\documentclass[twoside, 11pt]{article}

\usepackage{amssymb, amsmath, latexsym, mathrsfs, verbatim, calc, color}

\numberwithin{equation}{section}
\newtheorem{thm}{Theorem}[section]
\newtheorem{lem}[thm]{Lemma}
\newtheorem{cor}[thm]{Corollary}
\newtheorem{prop}[thm]{Proposition}
\newtheorem{rem}[thm]{Remark}

\newtheorem{defn}[thm]{Definition}

\def\ba{\begin{array}}
\def\ea{\end{array}}
\def\beq{\begin{equation}}
\def\endeq{\end{equation}}
\def\bes{\begin{equation*}}
\def\ees{\end{equation*}}
\def\bea{\begin{eqnarray}}
\def\eea{\end{eqnarray}}
\def\beaa{\begin{eqnarray*}}
\def\eeaa{\end{eqnarray*}}

\def\dis{\displaystyle}

\def\no{\noindent}

\def\lastline{\par \vspace{-7.3ex} \no}

\def\nts{\negthinspace}

\def\ms{\medskip}

\def\q{\quad}
\def\qq{\qquad}

\def\ol{\overline}

\def\={=\nts \nts=\nts \nts=\nts \nts=}

\def\({\textnormal{(}}
\def\){\textnormal{)}}

\def\cd{\cdot}
\def\cds{\cdots}

\def\pa{\partial}

\def\a{\alpha}

\def\d{\delta}
\def\e{\varepsilon}

\def\l{\lambda}
\def\m{\mu}
\def\n{\nu}

\def\si{\sigma}
\def\t{\tau}
\def\f{\varphi}
\def\th{\theta}
\def\o{\omega}

\def\f{\phi}
\def\vf{\varphi}

\def\D{\Delta}

\def\L{\Lambda}
\def\O{\Omega}

\def\cA{{\cal A}}

\def\cF{{\cal F}}

\def\cH{{\cal H}}

\def\cL{{\cal L}}

\def\cO{{\cal O}}

\def\cT{{\cal T}}

\def\dbF{\rm l\nts F}

\def\dbP{\rm l\nts P}

\def\hE{\mathbb{E}}
\def\hF{\mathbb{F}}

\def\hL{\mathbb{L}}

\def\hN{\mathbb{N}}

\def\hP{\mathbb{P}}
\def\hQ{\mathbb{Q}}
\def\hR{\mathbb{R}}

\def\sA{\mathscr{A}}

\def\sG{\mathscr{G}}

\def\sL{\mathscr{L}}
\def\sM{\mathscr{M}}

\def\esssup{\mathop{\rm esssup}}

\def\limsup{\mathop{\ol{\rm lim}}}

\def\qed{\hfill \rule[0cm]{.25cm}{.25cm}\medskip}   
\def\dfnn{\stackrel{\triangle}{=}}

\def\b1{{\bf 1}}
\def\neg{\negthinspace}
\def\dneg{\neg \neg}

\def\reff#1{{\rm(\ref{#1})}}

\newenvironment{itm}{\vspace{-1ex}\begin{itemize}}{\end{itemize}}
\def\bi{\begin{itm}}
\def\ei{\end{itm}}

\def\equ_ind{\arabic{section}.\arabic{equation}}
\def\sec_ind{\arabic{section}}

\begin{document}

\title{\textbf{Dynamic Equilibrium Limit
Order Book Model and Optimal Execution Problem}}

\author{Jin Ma\thanks{\noindent
 Department of Mathematics,
University of Southern California, Los Angeles, CA 90089. E-mail:
jinma@usc.edu. This author is supported in part by NSF grant
\#1106853. }, ~ Xinyang Wang\thanks{\noindent Institutional Equity Division, Morgan Stanley, New York, NY 10036. Email: xinyang.bryan.wang@gmail.com.}, ~
and ~{Jianfeng
Zhang}\thanks{\noindent Department of Mathematics, University of
Southern California, Los Angeles, CA 90089. E-mail:
jianfenz@usc.edu. This author is supported in part by NSF grant
\#1008873.}}

\date{}

\maketitle

\begin{abstract}  In this paper we propose a dynamic model of Limit Order Book (LOB).  The main feature of our model is
that the shape of the LOB is determined endogenously by
an expected utility function via a competitive equilibrium argument. Assuming zero resilience, the resulting
equilibrium density of the LOB is random, nonlinear, and time inhomogeneous. Consequently,
the liquidity cost can be defined dynamically in a natural way.

We next study an optimal execution problem in our model. We verify that the value function satisfies the Dynamic Programming Principle, and  is a viscosity solution to the corresponding
Hamilton-Jacobi-Bellman equation which is in the form of an integro-partial-differential quasi-variational inequality. We also prove the 
existence and analyze the structure of the optimal strategy via a verification theorem argument, assuming that the PDE has a 
classical solution. 
\end{abstract}

\vfill

\no{\bf Keywords:} Limit order book, liquidity risk, optimal execution, dynamic programming principle, viscosity solution, verification
theorem.

\no{\it 2000 AMS Mathematics subject classification:} 
91B51,70; 93E03, 20.

\section{Introduction}
\setcounter{equation}{0}

The effect of the liquidity of a security asset, both short term and long term, has been noticed by practitioners and researchers alike for quite some time. Tremendous efforts have been made in modeling the liquidity costs as well as its impact on the security prices
(see, e.g., \cite{Alfonsis,Alfonsis2,Cetins,OW}, to mention a few). In a frictionless market model
(Black-Scholes' framework for example), one assumes that
the securities can be bought or sold at a quote price regardless of the trade size and the
actual availability of the securities. But this is far from being realistic. In practice,
the parity between the supply and demand often causes the actual trade price to deviate from
the fundamental price, leading to the bid-ask spread. As a consequence, some extra cost
has to be paid in actual trading, especially when the volume of the trade is
relatively large compared to the existing liquidity on the market.

Unlike the quote driven market models, in which a market maker sets the price upon which all
the trades are made, an ``order-driven" market model is one that reflects more of the reality.
In such a model,
both buyers and sellers are allowed to be ``patient" in the sense that they submit the ``orders" containing the amount of the shares and the prices at which they are willing to buy or sell.
These orders are called {\it limit orders}. Unlike the ``market orders", which are executed immediately at the ``market price" whenever there is sufficient liquidity, the limit orders are executed only when an opposite order with the matching price comes in. Obviously, limit orders are usually not executed immediately, a {\it limit order book} (LOB) is thus formed. Intuitively, a reasonable
model of an LOB must contain the following basic elements:

(i) The best ask/bid price (the frontier of the sell/buy LOB);

(ii) The shape of the LOB (the volumes of the orders at each price).

There have been many papers in the literature trying to model and analyze the movement of the LOB (cf., e.g., \cite{Foucault,Handa,Hollifield,Rosu} and the references cited therein), as well as the optimal execution/liquidation problems in which a large trader needs to acquire/liquidate a certain amount of stocks in a given time horizon, with the minimal cost (see, e.g., \cite{AlfonsiSchied,Gatherals,OW}). Apart from the usual factors such as the
fundamental price (or mid-price) and the liquidity (often refer to the total amount of shares available for trading), an important characteristic  of an LOB is its ``shape", that is, the ``density" function of the LOB. This is particularly the case when the liquidity cost is among the main concerns. However, in most of the existing works the shape of the LOB is assumed to be {\it exogenously} given, either in the simple ``block-shaped" (cf. e.g., \cite{Alfonsis,OW}),  or in a general given shape that is supposed to be
determined by empirical studies (cf. e.g., \cite{AlfonsiSchied,Alfonsis2,shreve} and the references cited therein). However, such an assumption obviously lacks  the ability to adapt to the changes of market movement, especially when the underlying price is volatile within the concerned time horizon. A more ideal model would be such that the shape of the LOB could be determined {\it endogenously}, through some more basic market factors such as the bid-ask spread, the fundamental prices (the ``mid price", for example), and the market liquidity. This paper is an attempt in this direction.

To simplify the argument in this paper we shall consider
only the ``sell" side of the LOB, namely we assume that all the buyers are ``impatient" in
the sense that they only submit the market orders so there is no ``buy" side LOB. Our first objective is
to develop a dynamic model for the LOB whose shape is determined via the
movement of the fundamental price, the instantaneous trading size, as
well as the liquidity. The guiding principle of our model comes from  the idea of 
{\it equilibrium distribution}, initiated by Rosu \cite{Rosu}. Roughly speaking,
we assume  that there exists a competitive equilibrium among all the prices in
the LOB. The existence of such an equilibrium can be heuristically justified as a
balance between the expected sell price and the cost of waiting
(for the order to be executed). The equilibrium could be affected by the fundamental price, the execution of orders, and
the arrival of the new orders, etc., and when an existing equilibrium is broken, every seller in LOB will reposition
until an equilibrium is reached. It should be noted that this equilibrium is `` competitive" in the sense that one trader's deviation will be stopped by others' immediate undercutting. In other words, when the
market is under monopoly, we should allow the distribution to behave
differently.  In this paper we assume that the time of reaching new equilibrium is negligible, that is, the 
impact has zero duration, or ``zero resilience", for simplicity. We should note, however, that the issue of resilience
is interesting in its own right (see, e.g., \cite{Alfonsis} and also \cite{Alfonsis1, Alfonsis2}), but this is not the main purpose of this
paper.  

Mathematically, we shall assume that the equilibrium density process takes the form
$\mu^*_t=\mu^*(t,X_t,Q_t, y)$, $y\ge p_0$, where $p_0$ is the lowest (selling) price, $X$
is the fundamental value of the asset, and  $Q$ is the total volume of the LOB.
We also assume that the equilibrium is ``quantified' by a common expected utility on each price, which depends on the fundamental 
price and the total liquidity, and is denoted by $U(X, Q)$. Our main premise is that,  after each trade with size $\a \in [0, Q]$, the 
following two identities must hold:
\bea
\label{equil}
 \int_{p_0}^{p(\alpha)}\mu^{*}(X,Q,y)dy=\alpha, \q
 \frac{1}{\alpha}\int_{p_0}^{p(\alpha)}y {\mu}^{*}(X,Q,y)dy=U(X,Q-\alpha).
 \eea
Here the first equality is self-evident: $p(\alpha)=p(\alpha,X,Q)\ge p_0$ is the price in LOB at which
the accumulated volume of sell limit orders between $p_0$ and
$p(\a)$ is exactly equal to $\alpha$; whereas the second equality means that the
average price sold should be equal to $U(X,Q-\a)$, the expected utility for the remaining LOB
(a more detailed argument will be given in \S3).
Using the equations in (\ref{equil}) we will be able to solve explicitly the process $\mu^*$ in terms of $U$, and from which we will
define the liquidity cost, and argue that, modulo a term that is of
order $\a^2$, where $\a$ is the trading size, it is linear (although
time inhomogeneous) in $\a$. More importantly, once we obtain the density function of the LOB, we can
then evaluate the liquidity cost. We show that, under mild technical conditions, the average price (including liquidity cost)
exactly coincides with the {\it supply curve}  in sense of Cetin-Jarrow-Protter \cite{Cetins}.

Our second goal of this paper  is to consider an optimal execution problem. That is,
finding an optimal strategy of purchasing a large block of shares within a
prescribed time duration $[0,T]$ with a minimum cost. Such a problem has been studied by many authors (cf. e.g., 
\cite{Alfonsis2, Alfonsis, BayLud, shreve}, and the references cited therein), but with the endogenously given shape of LOB, our problem seems to be new. We shall consider only
two types of actions: the (buying) action of the large
investor self, and an aggregated action of all the other investors, which is modeled
as a compound Poisson process, representing all 
incoming limit sell orders, canceled orders, and the market buy orders. In other words, without the buying action of the investor, whose accumulated purchase will be described by an increasing process, the movement of the total available shares in the LOB is a
continuous time pure jump process.
We then show that the Bellman Principle of dynamic programming holds
in this case, and the value function is a viscosity solution of the
resulting HJB quasi-variational inequality (QVI). Finally, in the case that the QVI has a classical solution, we shall analyze the
optimal strategy by proving a verification theorem.  It is noted that 
the continuous (or inaction) region in our model may not be simply connected, and as a consequence the optimal strategy
may contain multiple (even infinitely many) jumps.

The rest of the paper is organized as follows. In \S \ref{sect-preliminary} we give the necessary technical background and describe the basic elements of the model. In \S \ref{sect-equilibrium} we introduce the notion of equilibrium distribution, and analyze some important quantities that can be derived endogenously from such distribution. These in particular include bid-ask spread and the liquidity cost that play the fundamental role in our optimal execution problem. In \S \ref{sect-execution} we introduce the optimal execution problem and study its various equivalent expressions.
In \S\ref{sect-DPP}  and \S \ref{sect-HJB} we prove the dynamic programming principle, derive the HJB equation, and prove that the value function is a viscosity solution to the corresponding HJB equation.  Finally, \S \ref{sect-strategy}  is devoted to the construction of
an optimal strategy, in the case that the HJB equation has a classical solution. 

\section{Preliminaries}
\setcounter{equation}{0}
\label{sect-preliminary}

Throughout this paper we assume that all the randomness comes from a complete
probability space $(\Omega, \mathcal{F}, \hP)$ on which are defined a standard Brownian motion
$W=\{W_{t}: t\ge 0\}$, and a standard Poisson process $N=\{N_{t}:\ge0\}$ with intensity $\l$. In what follows the Brownian motion $W$ represents the market noise that drives the fundamental value (or mid-price) of the underlying stock, and the Poisson process $N$ represents the frequency of the incoming limit orders. Therefore it is reasonable
to assume that $W$ and $N$ are independent. We shall denote
$\hF^{W}=\{\mathcal{F}^{W}_{t}:t\geq0\}$ and
$\hF^{N}=\{\mathcal{F}^{N}_{t}:t\geq0\}$ to be the natural
filtration generated by $W$ and $N$, respectively. Throughout the paper, we denote $\t_0:=0$ and let $0<\t_1<\t_2<\cds$ be the jump times of $N$.

We consider a finite time horizon $[0,T]$. For simplicity, we assume that there is only one stock traded in an order driven market, and the interest rate is 0.
We first give the mathematical description of the basic elements involved in our model.

\ms
\no{\bf 1. Fundamental Price.} We assume that the underlying stock has a {\it fundamental value} (or mid-price) which is known to the public. But the market price deviates away from it, due to the possible illiquidity, which leads to the bid-ask spread. Since the fundamental value only affects our
model as a source of randomness, we simply assume that it is a diffusion, and satisfies the following stochastic differential equation (SDE):
\bea
\label{X}
X_{t}=x+\int_0^tb(s,X_s)ds+\int_0^t\sigma(s,X_s)dW_s, \qq t\ge 0,
\eea
where $b$ and $\sigma$ satisfy the following standing assumptions: 

\begin{description}
\item[(H1)]  (i) $b(\cd, \cd)$ and $\si(\cd, \cd)$ are deterministic functions, continuous in $t$, and uniformly Lipschitz continuous in $x$, with a common uniform Lipschitz constant $L>0$.

(ii) $x>0$, $\si(t,0)=0$, and $b(t,0)\ge 0$. 
\end{description}

\begin{rem}
\label{rem-X}
{\rm It is clear that the assumption {\bf (H1)} guarantees the well-posedness of the the SDE (\ref{X}), and solution satisfies $X_t>0$ 
for all $t\ge 0$, $\hP$-a.s. The continuity of $b$ and $\si$ in $t$ is mainly for the viscosity property of the value function in \S \ref{sect-HJB} below.
For notational simplicity, in this paper we assume $W$ is $1$-dimensional, but all the results can be extended to higher dimensional case. Moreover, we may even allow $b$ and $\si$ to be random, and all the results in \S \ref{sect-execution} and \S \ref{sect-DPP} will still hold true, after obvious modification. However, in this case the HJB equation in Section \ref{sect-HJB} will become a backward stochastic PDE and the associated path dependent PDE. We refer to \cite{ETZ} for the related theory.
\qed}
\end{rem}

\no{\bf 2. The Limit Order Book (LOB).} We assume that there are patient and impatient
investors in the market, and they put different bid and/or ask prices to either liquidate or purchase the given stock based on their preferences (see \S3 for more discussion on this). Since in this paper we consider the optimal execution problem for purchasing the stock, only the sell side LOB will be relevant. We thus assume in what follows that all the buyers are impatient and only make ``market orders" (i.e., buying whatever is available on the market), and consequently there is
no ``buy side" LOB.  Moreover, we isolate one particular investor, referred as {\it the} investor, who will carry out the optimal execution problem later.

We shall assume that the movement of the LOB depends solely on the investment activities, namely the investor herself, and all other investors (buyers and sellers). For simplicity, we
assume that the activities of other investors are aggregated as a large investor whose
investment activities is described by a compound Poisson process $Y_t=\sum_{i=1}^{N_t}  {\L_i}$, $t\ge0$, where $\{\L_i\}_{i=1}^\infty$ is a sequence of i.i.d. random variables with distribution $\nu$. We shall assume $\hE\{|\L_i|\}<\infty$. 
We should note that the large investor is allowed to make both (buy and sell) limit orders and market orders, and can also cancel orders. Thus $\L_i$'s will take values in $\hR$ (i.e., $\D Y_t<0$ is possible). It is  useful to introduce the following filtration:
$\hF=\hF^{W}\otimes\hF^{Y}=\{\mathcal{F}^{W}_{t}\vee\mathcal{F}^{Y}_{t}:t\geq0\},$
which will be the basic information source allowed in our execution problem. We notice that $\hF^N\subset \hF^Y\subset \hF$.


\ms
\no{\bf 3. The Inventory Process.} We assume that the investor is trying to purchase a certain number, say $K$, shares of the given stock within a given time horizon $[0, T]$, and denote the accumulated number of shares up to time $t\in[0,T]$ by $\pi_t$. Then clearly
$\pi=\{\pi_t: t\ge0\}$ is an increasing process, and we assume that it is $\hF$-predictable. Note that, with this assumption, all the jumps times of $\pi$ is predictable, and consequently $\D \pi_{\t_i} \D Y_{\t_i}=0$, since all jump times of $N$ (and of $Y$) are totally inaccessible. In fact,
for practical reason we could, and will, assume that $N$ and $Y$ have c\`adl\`ag paths but  $\pi$ is c\`agl\`ad, and then naturally
we have
$$\D \pi_t \D Y_t := (\pi_{t+}-\pi_{t})(Y_{t}-Y_{t-})=0, \qq\forall t\in[0,T], \q
\mbox{$\hP$-a.s.}
$$
Note that with such a definition the investor can observe the jump of $Y$
and immediately jump afterwards. Clearly, each particular realization of $\pi$ could be considered as an execution strategy. We thus define
\bea
\label{calA}
\sA:=\{\pi: \mbox{~$\pi$ is $\hF$-predictable, non-decreasing, has c\`agl\`ad paths, and $\pi_T\le K$}\}.
\eea

\ms
We can now describe the dynamics of the total number of shares of the stock in the (sell) LOB, denoted by $Q=\{Q_{t}:t\in[0,T]\}$. We shall consider in this paper the
simplest case in which the dynamics of $Q$ can be affected by only two factors:
the order made by the investor herself, $\pi$, and the orders made by the other large investor (or the aggregated action by all other market participants), $Y$. Then, it is
readily seen that, for a given strategy $\pi\in\sA$ and initial inventory $q$, the movement of $Q^\pi := Q^{\pi, q}$  is determined by:
$Q^\pi_0 := q$, and 
\bea
\label{Qpi}
Q^\pi_t := Q^{\pi}_{\tau_{i}}-(\pi_{t}-\pi_{\tau_{i}})~\mbox{for}~ t\in(\tau_{i},\tau_{i+1}); \qq Q^{\pi}_{\tau_{i+1}} :=(Q^{\pi}_{\tau_{i+1}-}+\Delta
Y_{\t_{i+1}})^{+}.
\eea

\begin{rem}{\rm 
(i) $Q^\pi$ is c\`agl\`ad  in each $(\t_i, \t_{i+1})$. However, at $\t_{i+1}$, $Q^\pi$ can have left and/or right jumps. So $Q^\pi$ has both left and right limits, but in general it is neither left continuous nor right continuous on $[0, T]$.

(ii) When $\pi$ is continuous, which will be the  case in most of the paper, $Q^\pi$ is c\`adl\`ag.
\qed}
\end{rem}

We note from (\ref{Qpi}) that $Q^{\pi}_{\tau_{i+1}} \ge 0$. This is a natural constraint since the volume of the LOB can never be negative. However,  not all $\pi\in\sA$ will guarantee that the corresponding $Q^\pi_t \ge 0$ for all  $t\in [0, T]$. We thus consider the following {\it admissible strategies}: given $q\ge 0$,
\bea
\label{Aad}
\sA_{ad}(q) := \{\pi\in\sA: \mbox{~$Q^{\pi,q}_t\ge 0$ for all $t\in[0,T]$, $\hP$-a.s., where $Q^{\pi,q}$ is defined by \reff{Qpi}}\}.
\eea

Throughout the paper, we shall denote
\bea
\label{cO}
\hR_+ := (0,\infty),\q \bar \hR_+ := [0, \infty),\q \cO := \hR_+ \times [0, K) \times \hR_+,\q \bar\cO := \hR_+ \times [0, K] \times \bar\hR_+.
\eea
We remark that we do not take the closure for the first $\hR_+$ in $\bar\cO$.

\section{Equilibrium Distribution}
\label{sect-equilibrium} 

In this section we introduce the notion of ``equilibrium density" of the LOB, one of  the most important ingredients in our model. Our idea
follows from that of Rosu's \cite{Rosu}, which we now describe.
We assume that every seller comes into the market with the same amount of information
(this is different from the asymmetric information assumptions, cf. \cite{Avellaneda}). Each seller sets his/her ask price based on the personal preference, which
is the combination of the expected return of the order and the possible lost value (or cost) due to, say, the waiting time for the order to be executed. In an equilibrium we assume that every seller will have the same ``expected return" (or ``expected utility")
of the order, which we denote  by $U(X, Q)$, where $X$ is the fundamental value of the stock and $Q$ is the total number shares available.

The existence of such equilibrium could be argued as follows. Suppose two
sellers do not believe that they have the same expected return, then
one of them (usually the
one with lower expected return) is going to cancel his/her
limit order and resubmit it to the market with a different ask price in exchange
for a higher expected return. Then every seller in the market will do the same
until an equilibrium is reached. We should point out that such an equilibrium approach only works when there is sufficient competition in the market. In fact, when the
market is under monopoly, we should not expect the distribution to
behave like this.

Given the expected return $U(X,Q)$, we now introduce the concept of
``equilibrium density". Recall that the density function of an LOB
is a non-negative function $\mu(y)\ge 0$, $\forall y\ge 0$, such that
$\mu(y)=0$, for $y<p_0$, where $p_0\ge X$ is the lowest (best) ask
price, and that
\bea
\label{density}
\int_{p_0}^\infty \mu(y)dy=Q.
\eea

We note that if $\mu(y)\equiv \mu$, $p_0\le y\le p_0+Q\slash \mu$, is a constant, then the LOB is said to
have a ``block shape" (see, e.g., \cite{Alfonsis} and \cite{OW}).
Another way to study the problem is to assume the ``shape" of the LOB is given
exogenously (see, e.g., \cite{Alfonsis2, shreve}). Our main idea is to show
that the shape function is determined by the following simple facts.
Assume that a (large) market buy order comes in and $\a$-shares of the stock
were purchased, where $\a\in (0, Q]$. We assume that the lowest portion of $\alpha$
shares in the LOB is consumed. Thus, if we denote $p(0)=p(0,X,Q)$, to be
the lowest ask price, then we can find $p(\alpha)>p(0) $ such that
\bea
\label{pa}
\int_{p(0)}^{p(\alpha)}\mu(y)dy=\alpha.
\eea

On the other hand, we assume that, in equilibrium, 
 the average price of the sold block should have
the same expected return of the remaining  orders in the LOB, which has a total of 
$Q-\alpha$ shares after the purchase.
In other words, we assume that: for any $\alpha$ that $0\leq\alpha\leq Q$,
\bea
\label{dpa}
\frac{1}{\alpha}\int_{p(0)}^{p(\alpha)}y\mu(y)dy=U(X,Q-\alpha).
\eea

Now taking derivative with respect to $\alpha$ in (\ref{pa}) and (\ref{dpa}) we obtain:
\bea
\label{qp}
\left\{\ba{lll}
\mu(p(\alpha))p'(\alpha)=1;\ms\\
\dis \mu(p(\alpha))p'(\alpha)p(\alpha)=U(X,Q-\alpha)-\alpha\frac{\partial
U}{\partial x_{2}}(X,Q-\alpha).
\ea\right.
\eea
Solving two equations in (\ref{qp}) we have:
\bea
\label{palpha}
p(\alpha)&=&U(X,Q-\alpha)-\alpha\frac{\partial U}{\partial
Q}(X,Q-\alpha)
;\\
\label{qalpha}
\mu(p(\alpha))&=&\frac{1}{p'(\alpha)}=\left(\alpha\frac{\partial^{2}U}{\partial
Q^{2}}(X,Q-\alpha)-2\frac{\partial U}{\partial
Q}(X,Q-\alpha)\right)^{-1}.
\eea
We note that, by setting $\a=0$ in (\ref{palpha}),
\bea
\label{p0}
p(0, X, Q)=U(X,Q).
\eea
That is, the ``frontier" of the LOB is exactly the representative of the equilibrium, as
expected. On the other hand, since the function $\a\mapsto p(\a)$ is obviously non-decreasing, we can assume further that it is invertible and denote $h(y)=p^{-1}(y)$, then (\ref{qalpha}) becomes
\bea
\label{qalpha1}
\mu(y)&=&\frac{1}{p'(h(y))}=\left(h(y)\frac{\partial^{2}U}{\partial
Q^{2}}(X,Q-h(y))-2\frac{\partial U}{\partial
Q}(X,Q-h(y))\right)^{-1}.
\eea
Namely, the equilibrium density $\mu:= \mu^{X,Q}$ can be
explicitly derived, as long as $U(X,Q)$ is given.

We should remark here that the modeling of the expected return function $U(X, Q)$ is itself
an interesting and challenging problem. For example, in \cite{Rosu} such an expected return function was obtained explicitly by solving a recursive difference equation. Also, in
a slightly different setting, the relationship between the bid-ask spread and the
liquidity was considered by Avellaneda-Stoikov \cite{Avellaneda}, in which an argument
of indifference pricing was applied to construct the return function $U$.
In what follow we shall assume the existence of such a function $U$, and furthermore,
based on the discussion above, we make the following assumptions.

\begin{description}
\item{\bf (H2)} The expected utility function $U: \hR_+\times \bar \hR_+\mapsto \hR_+$ enjoys the following properties:

(i) $U$ is non-decreasing in $x$, and $\pa_Q U =\frac{\pa U}{\pa Q}< 0$, $\pa^2_Q U =\frac{\pa^2 U}{\pa Q^2}> 0$. 

(ii) $U$ is uniformly Lipschitz continuous in $(x,q)$, with  Lipschitz constant $L>0$.
%
\end{description}

\begin{rem}
{\rm (i) By \reff{qalpha}, the properties of $U$ in $q$  guarantees that $p'(\a)>0$, for all $0\le \a\le Q$, which leads further to the existence of its inverse so that the formula (\ref{qalpha1}) makes sense. Moreover,  by (\ref{p0}) we see that the function $p(0)=p(0,X,Q)$ is uniform Lipschitz for $(X, Q) \in \hR_+\times \bar\hR_+$. This fact will be frequently used in our discussion.

(ii) (H2) obviously does not render the function $U$ a true ``utility function" in either variable. In fact, the assumption (H2)-(i), 
which guarantees the positivity of the density function $\m$ (see (\ref{qalpha})), implies that it is decreasing and convex in $Q$,
hence a ``cost function" on $Q$  in a usual sense. Of course, it would be reasonable to 
 assume that $U$ is concave in $X$, hence a utility on the price, but we do not need such an assumption in the rest of our discussion.


(iii) In practice, it is natural to assume further that $U(x,q) \ge x$, or  $\lim_{q\to\infty}U(x,q)=x$.
 The latter implies that   the liquidity premium vanishes as the supply goes to infinity. But technically we do not need them in this paper.
\qed
}
\end{rem}

We conclude this section by observing that, given the density function $\mu=\mu^{X,Q}$, the cost for buying $\a$ shares of stock 
can be easily calculated as
\bea
\label{cost}
C(X,Q,\a):= \int_{p(0)}^{p(\a)}y\mu^{X,Q}(y)dy = \a U(X, Q-\a),
\eea
where the last equality is due to (\ref{dpa}). From this we obtain that
\bea
\label{liqcost}
\text{liquidity cost}&=&C(X,Q,\a)-\a X=[p(0)-X]\a+\int_{p(0)}^{p(\a)}[y-p(0)]\mu^{X,Q}(y)dy.
\eea
Clearly, we can see that the liquidity cost consists of a linear part (with respect to
the trade size $\a$), due to the bid-ask spread; and a higher order part that is determined by
the ``shape" of the LOB.
More precisely, assume for example $p'(\a)<\infty$, then we can easily derive from (\ref{liqcost}) that
\bea
\label{diffC}
C(X,Q,\a)=p(0)\a+O(\a^{2}).
\eea
In particular, if we consider a purchase strategy $\pi=\{\pi_t\}$, then (\ref{diffC}) amounts to saying that $C(X_t,Q^\pi_t,\D\pi_t)=p(0)\D\pi_t+O((\D\pi_t)^{2})$.
Consequently, for a continuous strategy $\pi^c=\{\pi^c_t, ~t\in[0,T]\}$, the following calculation
of the total cost will be useful in the rest of the paper:
\bea
\label{contcost}
\int_0^t C(X_s,Q^{\pi^c}_s,d\pi^c_s)=\int_0^tp(0,X_{s},Q_{s}^{\pi^c})d\pi^c_s = \int_0^tU(X_{s},Q_{s}^{\pi^c})d\pi^c_s .
\eea

\begin{rem}
{\rm The following obversion is worth noting. Assume that the function $U$ is sufficiently regular, then by (\ref{dpa})  we see that, for each $\a\in[0,Q]$, the process of ``average price" of the stock counting liquidity cost, defined by 
$$S(t,\o, \a)\dfnn \frac{1}{\a}C(X_{t}(\o),Q_{t}(\o),\a)=U(X_{t}(\o),Q_{t}(\o)-\a), \q (t,\o)\in [0,\infty)\times\O, $$
is a semi-martingale. Furthermore, the assumption (H2) implies that it is 
 convex and increasing with respect to the trade size $\a$. In other words,
the process $S$ is exactly the {\it supply curve} in the sense of Cetin-Jarrow-Protter \cite{Cetins}(!).
\qed}
\end{rem}

\section{Optimal Execution Problem}
\setcounter{equation}{0}
\label{sect-execution}

We are now ready to introduce the main objective of the paper: the optimal execution problem. Consider the
scenario when an investor would like to purchase  $K$ shares of the stock within
a prescribed time duration $[0,T]$.
Given initial inventory $q\ge 0$ and a purchase strategy $\pi\in\sA_{ad}(q)$, we consider the following cost functional:
\bea
\label{costfcnl}
J(\pi)=\mathbb{E}\Big\{\sum_{0\leq s<T}C(X_{s},Q^{\pi}_{s},\Delta\pi_{s})
+\int_{0}^{T}U(X_{s},Q^\pi_{s})d\pi^c_{s}+g(X_{T},K-\pi_{T})\Big\},
\eea
where $\pi^c$ denotes the continuous part of $\pi$, and $g: \hR_+ \times [0,K]\to \hR_+$ is the terminal penalty function. 
Clearly, the first term is the
cost for the jump part of $\pi$, and the
second term is the cost of the continuous part of $\pi$. The value function is thus
\bea
\label{valuefcn}
V_0:= V_0(q):=\inf_{\pi\in\sA_{ad}(q)}J(\pi).
\eea
We shall assume that the terminal penalty function $g$  satisfies the following assumption:
\begin{description}
\item[(H3)] (i) $g$ is uniformly Lipschitz continuous in $(x,y)$, with Lipschitz constant $L>0$.

(ii) For fixed $x$, $g$ is increasing and convex in $y$. Moreover, $g(x,0) =0$ and $g(x,y) \ge U(x,0) y$.
\end{description}

\begin{rem}
\label{rem-g}
{ \rm In the case $\pi_T < K$, one is forced to purchase the remaining amount of shares $y:=K-\pi_T$ at time $T$, regardless the 
liquidity.  The terminal (penalty)  $g(x,y) \ge U(x,0) y$ for $y\ge 0$ amounts to saying that this price would be more expensive than the highest market price $U(x,0)$, the price with zero liquidity. Furthermore, by (H3)-(ii) we see that $g(x,y) - g(x, y') \ge U(x,0) (y-y')$ for $0< y'<y$. Therefor if the final inventory is $Q$, and the investor needs to purchase a total of $y$ shares, but decides to buy $0< y'\le y\wedge Q$ from LOB right before $T$ and buys the remaining $y-y'$ using the penalty price, then his total cost would be: 
recall \reff{cost}, 
\beaa
C(x, Q, y') + g(x, y-y') = U(x,Q-y')y' + g(x, y-y') \le U(x,0) y' + g(x, y-y') \le g(x, y).
\eeaa 
This again shows that it is disadvantageous to purchase everything at the terminal time.
\qed}
\end{rem}

%

We now introduce two  alternative expressions for $V_0$ to facilitate the future discussion. First, we define the set of 
continuous strategies by
\bea
\label{sAc}
\sA^{c}_{ad}(q):=\{\pi\in\sA_{ad}(q): t\mapsto \pi_t ~\mbox{is continuous, $\hP$-a.s.}\}.
\eea
Clearly, if $\pi \in \sA^{c}_{ad}(q)$, then $Q^\pi$ is c\`adl\`ag and $C(X_t, Q^\pi_t, \D \pi_t) = 0$. We thus define
\bea
\label{V00}
\left\{\ba{lll}
\dis J^0(\pi) := \hE\Big\{\int_{0}^{T}U(X_{s},Q^\pi_{s})d\pi_{s}+g(X_{T},K-\pi_{T})\Big\}; \q \pi\in\sA_{ad}(q); \ms\\
\dis V^0_0 := \inf_{\pi\in\sA^c_{ad}(q)}J^0(\pi).
\ea\right.
\eea
Next, recall that $p(0, X, Q) = U(X, Q)$ is decreasing in $Q$. Thus, for $0<\a\le Q$, it holds that
\bea
\label{CgeD}
C(X,Q,\a)= \a U(X, Q-\a) = \int_0^\a U(X, Q-\a) du \ge \int_0^\a U(X,Q-u) du=: D(X, Q,\a). 
\eea
We now replace $C(\cds)$ by $D(\cds)$ in (\ref{costfcnl}) and define
\bea
\label{J1}
\left\{\ba{lll}
\dis J^1(\pi):=\mathbb{E}\Big\{\sum_{0\leq s<T}D(X_{s},Q^{\pi}_{s},\Delta\pi_{s})
+\int_{0}^{T}U(X_{s},Q^\pi_{s})d\pi^c_{s}+g(X_{T},K-\pi_{T})\Big\}, ~ \pi\in\sA_{ad}(q);\ms\\
\dis V^1_0:=\inf_{\pi\in\sA_{ad}(q)}J^1(\pi).
\ea\right.
\eea
%
We note that since $\sA^{c}_{ad}(q)\subseteq \sA_{ad}(q)$, it follows from   (\ref{CgeD})  that $V^1_0\le V_0
\le V^0_0$. Our main observation is that the cost $D(X, Q,\a)$ can actually be approximated by continuous strategies, thus these  inequalities should all be equalities. We substantiate this in the following theorem.
\begin{thm}
\label{thm-equiV0V1}
Assume (H1)- (H3). Then, it holds that $V^0_0=V_0=V^1_0$.
\end{thm}

{\it Proof.}
Since  $V_0^1\le V_0 \le V^0_0$ holds by definitions,  we need only show that $V^0_0\leq V^1_0$. To this end, we fix arbitrary $\pi\in \sA_{ad}(q)$ and $\e>0$. We claim that
\bea
\label{equi-claim}
 V^0_0 \leq J^1(\pi) + \e.
 \eea
 Indeed, for each $m\in\hN$, define $\t^m_0:=0$ and
$ \t^m_{i+1} := \inf\{t>\t^m_i: \D \pi_t \ge {1\over m}\}\wedge T$, $i=0,1,\cds.$
 Since $\pi$ has right limits and the filtration $\hF$ is right continuous, we see that $\t^m_i$ are $\hF$-stopping times, $\t^m_i < \t^m_{i+1}$ and $\D \pi_{\t^m_i} \ge {1\over m}$ whenever $\t^m_i < T$. Define
\bea
\label{pim}
\pi^{m}_s:= \pi^c_s + \sum_{i=1}^{m^2} \D \pi_{\t^m_i} 1_{\{\t^m_i \le s\}},  \qq s\in [0,T].
\eea
Clearly, $(\pi^m)^c=\pi^c$ and $\pi^m \le \pi$. This implies that $Q^{\pi^m} \ge Q^\pi$ and thus $\pi^{m}\in\sA_{ad}(q)$. Moreover, since $\sum_{i=1}^{m^2} \D \pi_{\t^m_i} \ge m$ on $\{\t^m_{m^2} < T\}$, we see that $\lim_{\m\to\infty}\hP(\t^m_{m^2} <T) = 0$. Consequently,
$\lim_{m\to\infty} \pi^m_T = \pi_T$, for all $\o$. Now by the monotonicity of $ U$ and \reff{CgeD}, we have
\beaa
\int_{0}^TU(X_{s},Q^{\pi^{m}}_s)d (\pi^m)^c_{s}&\leq&\int_{0}^TU(X_{s},Q^{\pi}_{s})d\pi^c_{s}; \\
\sum_{0\le s\le T}D(X_{s},Q^{\pi^{m}}_{s},\Delta\pi^{m}_{s}) &\le&\sum_{0\le s\le T}D(X_{s},Q^{\pi}_{s},\Delta\pi^{m}_{s}) = \sum_{i=1}^{m^2} D(X_{\t^m_i},Q^{\pi}_{\t^m_i},\Delta\pi_{\t^m_i}) \\
&\le& \sum_{0\le s\le T}D(X_{s},Q^{\pi}_{s},\Delta\pi_{s}).
\eeaa
Furthermore, since obviously one has $\lim_{m\to\infty} g(X_T, K-\pi^m_T) = g(X_T, K-\pi_T)$, we conclude that 
$\limsup_{m\to\infty} J^1(\pi^m) \le J^1(\pi)$, and thus there exists $M$ such that
\bea
\label{equi-est1}
J^1(\pi^M) \le J^1(\pi) +{\e\over 2}.
\eea

Next, recall again that $\D \pi_s \D N_s = 0$ and thus $\t_i \neq \t^M_j$, $\hP$-a.s. for all $i, j$. Let $\d>0$ be a small number.  For each $i=1,\cds, M^2$, let $j_i$ be the smallest $j$ such that $\t_j > \t^M_i$. We remark that $j_i$ is random and $\t_{j_i}$ is still an $\hF$-stopping time. Define $\pi^{M,\d}$ recursively as follows. First, $\pi^{M,\d}_s := \pi^c_s$ for $0\le s\le \t^M_1$. For $i=1,\cds, M^2$, denote $\t^{M,\d}_i := [\t^M_i+\d] \wedge \t^M_{i+1} \wedge \t_{j_i}$ , and define
\bea
\label{pimd}
\pi^{M, \d}_s := \pi^{M,\d}_{\t^m_i} + [\pi^c_s - \pi^c_{\t^M_i}] + {s\wedge \t^{M,\d}_i-\t^M_i\over \d} \D \pi_{\t^M_i}, \q s\in (\t^M_i, \t^M_{i+1}], 
\eea
where we abuse the notation that $\t^M_{m^2+1} := T$. It is clear that $\pi^{M,\d}$ is continuous and $\pi^{M,\d} \le \pi^M$. This implies that $\pi^{M,\d} \in \sA^c_{ad}(q)$. Note that, by changing variable $u := \t^M_i + {\a\over \d}(s-\t^M_i)$, we have
$$D(X, Q,\a)  = \int_0^\a U( X, Q-u) du ={\a\over \d} \int_{\t^M_i}^{\t^M_i+\d} U(X, Q- {\a\over \d} (s-\t^M_i)) ds.
$$
On the other hand, it is not hard to check that, for $s\in [\t^M_i, \t^{M,\d}_i]$, it holds that
\beaa
&&Q^{\pi^{M, \d}}_s =  Q^{\pi^{M,\d}}_{\t^M_i} -   [\pi^c_s - \pi^c_{\t^M_i}] - {s-\t^M_i\over \d} \D \pi_{\t^M_i} \ \ge  Q^{\pi}_{\t^M_i} -   [\pi^c_s - \pi^c_{\t^M_i}] - {s-\t^M_i\over \d} \D \pi_{\t^M_i}, 
\eeaa
and that $\lim_{\d\to 0} \hP(\t^{M,\d}_i = \t^M_i + \d) = 1$, thus  we have $\lim_{\d\to 0} \pi^{M,\d}_T = \pi_T$, $\hP$-a.s.

Now, by the monotonicity of $U$ again and applying the dominated convergence theorm,
\beaa
&&J^0(\pi^{M,\d}) -  J^1(\pi^M) \\
&=& \hE\Big\{ \int_0^T [U(X_s, Q^{\pi^{M,\d}}_s) - U(X_s, Q^\pi_s)]d \pi^c_s + [g(X_T, K-\pi^{M,\d}_T) - g(X_T, K-\pi_T)]\\
&&+ \sum_{i=1}^{M^2} [\int_{\t^M_i}^{\t^{M,\d}_i} {\D \pi_{\t^M_i}\over \d} U( X_s, Q^{\pi^{M,\d}}_s) ds - D(X_{\t^M_i}, Q^\pi_{\t^M_i}, \D \pi_{\t^M_i})]\Big\} \\
&\le&\hE\Big\{ [g(X_T, K-\pi^{M,\d}_T) - g(X_T, K-\pi_T)]+ \sum_{i=1}^{M^2}   \int_{\t^M_i}^{\t^{M}_i+\d} {\D \pi_{\t^M_i}\over \d}\times\\
&&\Big[p\big(0, X_s, Q^{\pi}_{\t^M_i} - (\pi^c_s-\pi^c_{\t^M_i}) - {\D \pi_{\t^M_i}\over \d} (s-\t^M_i)\big) -p\big(0, X_{\t^M_i}, Q^{\pi}_{\t^M_i} - {\D \pi_{\t^M_i}\over \d} (s-\t^M_i)\big)\Big] ds\Big\}\\
&\le& L\hE\Big\{ |\pi^{M,\d}_T-\pi_T|+ {\pi_T\over \d} \sum_{i=1}^{M^2}   \int_{\t^M_i}^{\t^{M}_i+\d} \Big[|X_s - X_{\t^M_i}| +| \pi^c_s-\pi^c_{\t^M_i}|\Big] ds\Big\} \to 0, ~\mbox{as}~\d\to 0.
\eeaa
Setting $\d>0$ small enough such that  $J^0(\pi^{M,\d}) \le  J^1(\pi^M) + {\e\over 2}$. By \reff{equi-est1} and recalling that $\pi^{M,\d} \in \sA^c_{ad}(q)$, we prove \reff{equi-claim}, whence the theorem.
\qed

We conclude this section with a dynamic version of the value function $V$.  Let $(t,x,k,q)\in [0,T]\times \bar\cO$  (recall \reff{cO}), 
and let $X^{t,x}$ be the solution to  \reff{X} on $[t, T]$ with $X_t = x$, a.s. Denote
$$\sA(t,k):=\{\pi: \mbox{$\pi$ is $\hF$-predictable, c\`agl\`ad, non-decreasing, $\pi_t =k$,  and $\pi_T\le K$}\}. 
$$
Denote $\t^t_0 := t$, and $\t^t_i$, $i\ge 1$, being the jump times of $N$ on $(t, T]$.  For any $\pi \in \sA(t,k)$,  let 
\bea
\label{Qpit}
\left\{\ba{lll}
Q^\pi_s := Q^{\pi}_{\tau^t_{i}}-(\pi_{s}-\pi_{\tau^t_{i}}) \qq \mbox{for}~ s\in(\tau^t_{i},\tau^t_{i+1});\\
Q^\pi_t := q;\q Q^{\pi}_{\tau^t_{i+1}} :=(Q^{\pi}_{\tau^t_{i+1}-}+\Delta Y_{\t^t_{i+1}})^{+}, \qq i\ge 1,
\ea\right.
\eea
and define 
\bea
\label{sAt}
\left.\ba{lll}
\sA_{ad}(t,k,q) := \{\pi \in \sA(t,k): Q^{\pi, q}_s \ge 0,~s\in [t, T], \hP\mbox{-a.s.}\}, \\
\sA^c_{ad}(t,k,q) := \{\pi \in \sA_{ad}(t,x,q): \pi ~\mbox{is continuous, $\hP$-a.s.}\}.
\ea\right.
\eea
By Theorem \ref{thm-equiV0V1}, we now define the dynamic value function $V$ via two equivalent expressions:
\bea
\label{Vt}
V(t,x,k,q) &:=&\inf_{\pi\in\sA^c_{ad}(t,k,q)}J^0(t,x,k,q;\pi) = \inf_{\pi\in\sA_{ad}(t,k,q)}J^1(t,x,k,q;\pi),
\eea
where
\bea
\label{J0J1}
J^0(t,x,k,q;\pi)&:=& \hE\Big\{\int_{t}^{T}U(X^{t,x}_{s},Q^\pi_{s})d\pi_{s}+g(X^{t,x}_{T},K-\pi_{T})\Big\};\\
J^1(t,x,k,q;\pi)&:=&\mathbb{E}\Big\{\sum_{t\leq s<T}D(X^{t,x}_{s},Q^{\pi}_{s},\Delta\pi_{s})
+\int_{t}^{T}U(X^{t,x}_{s},Q^\pi_{s})d\pi^c_{s}+g(X^{t,x}_{T},K-\pi_{T})\Big\}.\nonumber
\eea

\begin{rem}
\label{rem-J01}
 {\rm (i) We note that the cost functional $J^0(t,x,k,q;\pi)$ in (\ref{J0J1}) uses only continuous strategies. It will facilitate the argument when
we prove that the value function $V$ is a viscosity solution to the HJB equation in \S \ref{sect-DPP} and \S \ref{sect-HJB}. 

(ii) The cost functional  $J^1(t,x,k,q;\pi)$ will be useful when we investigate the existence of optimal strategy in \S \ref{sect-strategy}. 
Recall from Theorem \ref{thm-equiV0V1} the inequality
$V^0_0\le V_0 \le V^1_0$. Thus 
an optimal strategy, if exists, should also optimize $J^1$. 
However, it is worth noting that cost function $D(\cds)$ does not have a practical meaning,  
as opposed to the cost function $C(\cds)$, and in practice it  cannot be implemented directly.
Nevertheless, combining the approximations \reff{pim} and \reff{pimd} in  the proof of Theorem \ref{thm-equiV0V1}, we will be
able to find an implementable good approximation of optimal strategy, as we shall see in  \S \ref{sect-strategy}. 
\qed}
\end{rem}

\section{Dynamic Programming Principle}
\setcounter{equation}{0}
\label{sect-DPP}

In this section we verify some properties of the value function $V$ and establish the Dynamic Programming Principle (DPP).  As we pointed out in Remark \ref{rem-J01}-(i), we shall consider the cost functional $J^0$. 
We begin by the  regularity of $V$ with respect to the ``spatial variables" $x$, $k$, and $q$, respectively.
\begin{prop}
\label{prop-Vregxkq}
Assume (H1)-(H3). Then for each $t\in [0,T]$, the value function $V(t,x,k,q)$ is 
non-decreasing $x$, non-increasing in $k$ and $q$, respectively, and uniformly Lipschitz continuous 
with respect to $(x,k,q)\in \bar \cO$.
\end{prop}

{\it Proof.} We first check the properties in $x$. Assume $x_{1}< x_{2}$. Then by the comparison theorem of SDE, we have 
$X^{t,x_1}_s \le X^{t,x_2}_s$,  for all $t\le s\le T$, $\hP$-a.s. Since both $U$ and $g$ are non-decreasing and uniformly Lipschitz continuous in $x$,  for any $\pi\in \sA^c_{ad}(t,k,q)$ we see that 
\bea
\label{J0x}
0&\le&J^0(t,x_{2},k,q;\pi)-J^0(t,x_{1},k,q;\pi)\\
&=&\mathbb{E}\Bigg\{\int_{t}^{T}[U(X^{t,x_2}_{s},Q^{\pi}_{s})-U(X^{t,x_1}_{s},Q^{\pi}_{s})]d\pi^c_{s}
+g(X^{t,x_{2}}_{T},K-\pi_T )-g(X^{t,x_{1}}_{T},K-\pi_T)\Bigg\}\nonumber\\
&\leq&C\mathbb{E}\left\{\max_{s\in[t,T]}|X^{t,x_{2}}_{s}-X^{t,x_{1}}_{s}|\right\}\le C (x_{2}-x_{1}).\nonumber
\eea
Switching the role of $x_1$ and $x_2$ we can easily deduce the Lipschitz property in $x$:
\bea
\label{Lipx}
|V(t,x_{2},k,q)-V(t,x_{1},k,q)|\leq C|x_{2}-x_{1}|, \qq \forall x_1, x_2\in \hR.
\eea
%

We next check the properties in $k$. Let $0\le k_{1} < k_{2}\le K$. For any $\pi\in\sA^c_{ad}(t,k_1,q)$, consider
the strategy $\pi'_s:=[k_2+(\pi_s-k_1)]\wedge K$, $s\in[t,T]$. 
Clearly, $\pi'\in\sA^c_{ad}(t,k_2,q)$, and it satisfies: $\pi'_T \ge \pi_T$,  $d \pi'_s \le d\pi_s$ $s\in[t, T]$. Consequently
we have $Q^{\pi', q} \ge Q^{\pi,q}$,  $J^0(t,x,k_2, q;\pi')\leq J^0(t,x,k_1,q;\pi)$, and thus  $V(t, x, k_{2},q)
\leq V(t,x, k_{1},q)$.
On the other hand, for any strategy
$\pi\in\sA^c_{ad}(t,k_{2},q)$, let $\pi':=\pi-(k_2-k_1)\in\sA^c_{ad}(t,k_{1},q)$. Then  $Q^{\pi', q} = Q^{\pi,q}$, and thus:
\bea
\label{J0k}
J^0(t,x,k_{1},q;\pi')-J^0(t,x,k_{2},q;\pi)
=\mathbb{E}\left\{g(X^{t,x}_{T},K-\pi'_T)-g(X^{t,x}_{T},K-\pi_T)\right\}\le C(k_{2}-k_{1}).
\eea
Similar to (\ref{Lipx} this implies  the uniform Lipschitz continuity of $V$ in $k$.
%

It remains to prove the Lipschitz property in $q$.  As before we first assume $0\le q_1 < q_2$. It is clear that $\sA^c_{ad}(t,k,q_{1}) \subset \sA^c_{ad}(t,k,q_{2})$, and for any $\pi \in \sA^c_{ad}(t,k,q_{1})$, we have $Q^{\pi,q_1}_s \le Q^{\pi, q_2}_s$. Then 
\bea
\label{J0q}
J^0(t,x,k, q_1;\pi) \ge J^0(t,x,k, q_2;\pi)  &\mbox{for all}& \pi \in \sA^c_{ad}(t,k,q_{1}),
\eea
which leads to   $V(t,x,k, q_1) \ge V(t,x,k, q_2)$. On the other hand, note that $\pi^0\equiv 
 k \in  \sA^c_{ad}(t,k,q_{1})$. For any $\pi \in \sA^c_{ad}(t,k,q_{2})$, denote $\D Q := Q^{\pi, q_2}- Q^{k,q_1}$ and $\t := \inf\{s\ge t: \D Q_s \le 0\} \wedge T$.  Recall \reff{Qpit},   by induction on $i$ one deduce easily  that $\D Q$ is non-increasing on $[t,\t]$. Then
\bea
\label{Dpi}
\pi_\t - \pi_t &=& \sum_{i=0}^\infty [\pi_{\t \wedge \t^t_{i+1}} - \pi_{\t^t_i}]1_{\{\t^t_i <\t\}} = \sum_{i=0}^\infty [ \D Q_{ \t^t_i} - \D Q_{\t\wedge \t^t_{i+1}-}]1_{\{\t^t_i <\t\}}\nonumber\\
& \le& \D Q_t - \D Q_{\t-} \le \D Q_t = q_2-q_1. 
\eea
Now define $\pi'_s := \pi^0_s 1_{[t,\t]}(s) + [\pi_s - \pi_\t]1_{(\t, T]}$. Since $\pi$ is continuous and $\pi^0\equiv k$, by \reff{Qpit} we see that $\D Q_\t = 0$, as $\t<T$. Then $Q^{\pi', q_1}_s = Q^{k, q_1}_s \le Q^{\pi, q_2}_s$, $s\in [t, \t]$, and $Q^{\pi', q_1}_s = Q^{\pi, q_2}_s$, $s\in (\t, T]$. Namely $\pi' \in \sA^c_{ad}(t,k,q_{1})$.
Moreover, \reff{Dpi} implies that $0 \le \pi_T - \pi'_T = \pi_\t-\pi_t \le q_2-q_1$. Then
\beaa
&&J^0(t,x,k, q_1;\pi') - J^0(t,x,k, q_2;\pi)\\
&=&\mathbb{E}\Big\{-\int_{t}^{\t}U(X^{t,x}_{s},Q^{\pi,q_2}_{s})d\pi_{s} +g(X^{t,x}_{T},K-\pi'_T )-g(X^{t,x}_{T},K-\pi_T)\Big\}\\
&\leq&C\mathbb{E}\left\{\pi_T-\pi'_T\right\}\le C (q_{2}-q_{1}).
\eeaa
Since $\pi \in \sA^c_{ad}(t,k,q_{2})$ is arbitrary, we obtain $V(t,x,k, q_1) - V(t,x,k, q_2)\le C (q_{2}-q_{1})$. Reversing the role
of $q_1$ and $q_2$ we obtain the Lipschitz property of $V$ in $q$, proving the proposition.
\qed

We can now follow the standard arguments in the literature to establish the following simpler from of dynamic programming principle, when the time increments are deterministic.
\begin{prop}
\label{prop-DPPProp}
Assume (H1) - (H3). Then, for any $0\le t_1 < t_2 \le T$ and $(x,k,q)\in\bar\cO$, 
\bea
\label{DPP0}
V(t_1,x,k,q)=\inf_{\pi\in\sA^c_{ad}(t_1,k,q)}
\mathbb{E}\left\{\int_{t_1}^{t_2}U(X^{t_1,x}_{s},Q^{\pi,q}_{s})d\pi_{s}+V(t_2,X^{t_1,x}_{t_2},\pi_{t_2},Q^{\pi,q}_{t_2})\right\}.
\eea
\end{prop}

{\it Proof.} Let $\tilde V(t_1,x,k,q)$ denote the right side of \reff{DPP0}. 
We first show that $ V(t_1,x,k,q)\ge \tilde V(t_1,x,k,q)$. Indeed, for any $\pi\in\sA_{ad}^c(t_1,k,q)$, let $\tilde\pi$ denote the restriction of $\pi$ on $[t_2, T]$. Then 
$X^{t_2,X^{t_1,x}_{t_2}}_s = X^{t_1,x}_s$, $Q^{\tilde\pi, Q^{\pi,q}_{t_2}}_s = Q^{\pi,q}_s$, for $s\in [t_2, T]$.
In other words, $ \tilde\pi \in \sA_{ad}^c(t_2, \pi_{t_2}, Q^{\pi, q}_{t_2})$. This implies that
\beaa
&&J^0(t_1,x,k,q;\pi)=\hE\Big\{\int_{t_1}^{T}U(X^{t_1,x}_{s},Q^{\pi,q}_{s})d\pi_{s}+g(X^{t_1,x}_{T},K-\pi_{T})\Big\}\\
&=&\hE\Big\{\int_{t_1}^{t_2}U(X^{t_1,x}_{s},Q^{\pi,q}_{s})d\pi_{s}+\hE\Big[\int_{t_2}^{T}U(X^{t_2,X^{t_1,x}_{t_2}},Q^{\tilde\pi, Q^{\pi,q}_{t_2}}_{s})d\pi_{s}+g(X^{t_2,X^{t_1,x}_{t_2}}_{T},K-\tilde\pi_{T})\Big|\cF_{t_2}\Big]\Big\}\\
&=&\hE\Big\{\int_{t_1}^{t_2}U(X^{t_1,x}_{s},Q^{\pi,q}_{s})d\pi_{s}+J^0(t_2,X^{t_1,x}_{t_2},\pi_{t_2},Q^{\pi,q}_{t_2}; \tilde\pi)\Big\}\\
&\ge& \hE\Big\{\int_{t_1}^{t_2}U(X^{t_1,x}_{s},Q^{\pi,q}_{s})d\pi_{s}+V(t_2,X^{t_1,x}_{t_2},\pi_{t_2},Q^{\pi,q}_{t_2})\Big\}.
\eeaa
We remark that in the above the last equality can be proved rigorously by using the notion of regular conditional probability distribution. Since the argument would be rather lengthy but more or less standard, we omit the details.  Now take infimum  over $\pi\in\sA_{ad}^c(t_1,k,q)$ on both sides of above, we obtain $ V(t_1,x,k,q)\ge \tilde V(t_1,x,k,q)$.

To prove the opposite inequality, we first fix $\e>0$, and consider  a countable partition $\{O_{i}\}_{i=1}^{\infty}$  of
$\bar\cO$  and $(x_i,k_i,q_i)\in O_i$, $i=1,2 \cds$,  such that, for any $(x,k,q)\in O_i$, it holds that $|x-x_i|\le \e$, $k_i-\e\leq k \le k_i$, and $q_i\leq q\le q_i+\e$.  Now for each $i$,  choose $\pi^{i}\in\sA^c_{ad}(t_2,k_{i},q_{i})$ such that
\beaa
J^0(t_2,x_{i},k_{i},q_{i};\pi^{i})\le V(t_2,x_{i},k_{i},q_{i}) +\e.
\eeaa
For any $(x,k,q)\in O_{i}$, note that $\pi^i - k_i + k \in \sA^c_{ad}(t_2,k,q_{i}) \subset \sA^c_{ad}(t_2,k,q)$. Then, by \reff{J0x}, \reff{J0k}, \reff{J0q}, and applying  Proposition \ref{prop-Vregxkq}, for a generic constant $C$ we have
\bea
\label{J0eps}
J^0(t_2,x,k,q;\pi^{i}-k_i+k)&\leq&J^0(t_2,x_{i},k_i,q;\pi^{i})+C\e \le J^0(t_2,x_{i},k_i,q_i;\pi^{i})+C\e  \nonumber\\
&\leq&V(t_2,x_{i},k_i,q_i)+C\e \le V(t_2,x,k,q)+C\e.
\eea
Now for any $\pi\in\sA^c_{ad}(t_1,k,q)$, define a new strategy $\tilde{\pi}$:
\beaa
\tilde \pi_s := \pi_s 1_{[t_1, t_2]}(s) +  \Big[\sum_i [\pi^i_s - k_i + \pi_{t_2}] 1_{D_i}(X^{t_1,x}_{t_2}, \pi_{t_2}, Q^{\pi,q}_{t_2})\Big] 1_{(t_2, T]}(s).
\eeaa
It is clear that $\tilde\pi_{t_1} = k$, $\tilde\pi$ is continuous and non-decreasing on $[t, T]$,  and $\tilde\pi_T  \le \pi^i_T \le K$ on each $O_i$. Moreover, $Q^{\tilde\pi, q}_s = Q^{\pi,q}_s \ge 0$ for $s\in [t_1, t_2]$, and for  $s\in [t_2, T]$, on $O_i$ we have
\beaa
Q^{\tilde\pi, q}_s = Q^{\pi^i, Q^{\pi,q}_{t_2}}_s \ge Q^{\pi^i, q_i}_s \ge 0.
\eeaa
Thus  $\tilde{\pi}\in\sA^c_{ad}(t_1,k,q)$, and therefore, it follows from  \reff{J0eps} that
\beaa
&&V(t_1,x,k,q)\leq J^0(t_1,x,k,q;\tilde{\pi})\\
&=&\hE\Big\{\int_{t_1}^{t_2}U(X^{t_1,x}_{s},Q^{\pi,q}_{s})d\pi_{s} +\hE\Big[\int_{t_2}^{T}U(X^{t_1,x}_s,Q^{\tilde{\pi},q}_{s})d \tilde \pi_{s}
+g(X^{t_1,x}_{T},K-\tilde{\pi}_{T})\Big|\cF_{t_2}\Big]\Big\}\\
&=&\hE\Big\{\int_{t_1}^{t_2}U(X^{t_1,x}_{s},Q^{\pi,q}_{s})d\pi_{s} +J^0(t_2,X^{t_1,x}_{t_2},\pi_{t_2},Q^{\pi,q}_{t_2};\tilde{\pi})\Big\}\\
&=&\hE\Big\{\int_{t_1}^{t_2}U(X^{t_1,x}_{s},Q^{\pi,q}_{s})d\pi_{s} + \sum_i J^0(t_2,X^{t_1,x}_{t_2},\pi_{t_2},Q^{\pi,q}_{t_2};\pi^i-k_i + \pi_{t_2})1_{D_i}(X^{t_1,x}_{t_2}, \pi_{t_2}, Q^{\pi,q}_{t_2})\Big\}\\
&\le &\hE\Big\{\int_{t_1}^{t_2}U(X^{t_1,x}_{s},Q^{\pi,q}_{s})d\pi_{s} + V(t_2,X^{t_1,x}_{t_2},\pi_{t_2},Q^{\pi,q}_{t_2}) \Big\} + C\e,
\eeaa
Now, since  $\e>0$ is arbitrary and $\pi\in\sA^c_{ad}(t_1,k,q)$, we conclude  that $ V(t_1,x,k,q)\le \tilde V(t_1,x,k,q)$, proving the proposition.
\qed

As a corollary of Proposition \ref{prop-DPPProp}, we shall prove the temporal regularity of $V$. We note that this will be a crucial 
step towards the general form of dynamical programming principle.
\begin{cor}
\label{cor-Vregt}
Assume (H1)-(H3). Then,  for any $0\le t_1 <t_2\le T$ and $(x,k,q) \in\bar\cO$, we have
\bea
\label{Vreg}
|V(t_1,x,k,q) - V(t_2,x,q)|\le C(1+|x|) \sqrt{t_2-t_1}.
\eea
\end{cor}
 {\it Proof.} First note that the constant process $k\in \sA_{ad}^c(t_1,k,q)$.  Then, by Propositions \ref{prop-DPPProp} and \ref{prop-Vregxkq}, 
 \beaa
V(t_1,x,k,q) - V(t_2,x,k,q)&\le& \hE \{V(t_2,X^{t_1,x}_{t_2}, k,Q^{k,q}_{t_2}) \}- V(t_2,x,k,q)\\
&\le& C \mathbb{E}\{|X^{t_1,x}_{t_2}-x|+|Q^{k,q}_{t_2}-q|\}.
\eeaa
Next, recall from \S2 that the dynamics of $Q$ (see (\ref{Qpi})) is driven by the compound Poisson process $Y$, whose jump size $\L_i$'s  and the
jump times $\t_i$'s  are independent. Then one can easily check:
\bea
\label{Vregt1}
\mathbb{E}\{|X^{t_1,x}_{t_2}-x|\}&=&  \mathbb{E}\Big\{ \Big|\int_{t_1}^{t_2} b(s, X^{t_1,x}_s) ds + \int_{t_1}^{t_2} \si(s, X^{t_1,x}_s) dW_s\Big|\Big\}\le C(1+|x|) \sqrt{t_2-t_1};\nonumber\\
\mathbb{E}\{|Q^{k,q}_{t_2}-q|\} &\le&\mathbb{E}\Big\{ \sum_{i=1}^\infty |\L_i| 1_{\{t_1< \t_i \le t_2\}}\Big\} =\sum_{i=1}^\infty \mathbb{E}\{  |\L_i| \} \mathbb{E}\{1_{\{t_1< \t_i \le t_2\}}\}  \\
&=& \hE\{|\L_1|\} \mathbb{E}\Big\{ \sum_{i=1}^\infty  1_{\{t_1< \t_i \le t_2\}}\Big\}=  \hE\{|\L_1|\} \mathbb{E}\Big\{ N_{t_2} - N_{t_1}\Big\} = \l  \hE\{|\L_1|\} [t_2-t_1].\nonumber
\eea
Consequently, we obtain 
 \bea
 \label{Vregt2}
V(t_1,x,k,q) - V(t_2,x,k,q) \le C(1+|x|) \sqrt{t_2-t_1}.
 \eea
On the other hand, since $U\ge 0$ and $V$ is decreasing in $q$,
 \beaa
V(t_2,x,k,q) - V(t_1,x,q)
&\leq&V(t_2,x,k,q) - \inf_{\pi\in\sA^c_{ad}(t_1,k,q)} \mathbb{E}\left\{V(t_2,X^{t_1,x}_{t_2},\pi_{t_2},Q^{\pi,q}_{t_2})\right\}\\
&=& \sup_{\pi\in\sA^c_{ad}(t_1,k,q)} \mathbb{E}\left\{V(t_2,x,k,q) - V(t_2,X^{t_1,x}_{t_2},\pi_{t_2},Q^{\pi,q}_{t_2})\right\}\\
&\le& C  \sup_{\pi\in\sA^c_{ad}(t_1,k,q)}\mathbb{E}\Big\{|X^{t_1,x}_{t_2}-x|+[Q^{\pi,q}_{t_2}-q]^+\Big\}\\
&=& C \mathbb{E}\Big\{|X^{t_1,x}_{t_2}-x|+[Q^{k,q}_{t_2}-q]^+\Big\}
\le C(1+|x|) \sqrt{t_2-t_1},
 \eeaa
 where the last inequality is due to \reff{Vregt1}. This, together with \reff{Vregt2}, leads to \reff{Vreg}.
\qed

To conclude this section we give a general version of the dynamic programming principle. Denote $\cT_t$ to be 
all the $\hF$-stopping times taking values in $(t,T]$. 
\begin{thm}
\label{thm-DPP}
Assume (H1)-(H3). Then, for any $(t,x,k,q)\in [0,T)\times \bar\cO$ and any  $\t\in \cT_t$, 
\bea
\label{DPP}
V(t,x,k,q)=\inf_{\pi\in\sA^c_{ad}(t,k,q)}
\mathbb{E}\left\{\int_{t}^{\tau}U(X^{t,x}_{s},Q^{\pi,q}_{s})d\pi_{s}+V(\tau,X^{t,x}_{\tau},\pi_{\tau},Q^{\pi,q}_{\tau})\right\}.
\eea
\end{thm}
{\it Proof.} For each $\pi\in\sA^c_{ad}(t,k,q)$ and $\t\in \cT_t$, denote $I(\pi,\t)$ be the expectation on the right side of \reff{DPP}.  Following the arguments in Proposition \ref{prop-DPPProp} one can easily show that $ V(t,x,k,q)\ge \inf_{\pi\in\sA^c_{ad}(t,k,q)} I(\pi,\t)$. So it suffices to prove the reversed inequality: 
\bea
\label{DPPineq}
V(t,x,k,q)\le \inf_{\pi\in\sA^c_{ad}(t,k,q)} I(\pi,\t).
\eea

We first assume that $\t\in\cT_t$ takes only finitely many values $t< t_1<\cds <t_m\le T$. We prove \reff{DPPineq} by induction on $m$. When $m=1$, \reff{DPPineq} follows from Proposition \ref{prop-DPPProp}. Now assume that \reff{DPPineq} holds for $m-1$, and
that $\t$ takes $m$ values. For any $\pi \in\sA^c_{ad}(t,k,q)$, we have  
\bea
\label{Ipit}
I(\pi,\t)
&=& \mathbb{E}\Big\{\int_{t}^{t_1}U(X^{t,x}_{s},Q^{\pi,q}_{s})d\pi_{s}+V(t_1,X^{t,x}_{t_1},\pi_{t_1},Q^{\pi,q}_{t_1}) 1_{\{\t=t_1\}} \nonumber\\
&&+ \Big[\int_{t_1}^{\tau}U(X^{t,x}_{s},Q^{\pi,q}_{s})d\pi_{s}+V(\tau,X^{t,x}_{\tau},\pi_{\tau},Q^{\pi,q}_{\tau})\Big]1_{\{\t > t_1\}} \Big\}. \nonumber 
\eea
Note that $\{\t>t_1\} \in \cF_{t_1}$ and $\t$ takes only $m-1$ values on $\{\t>t_1\}$. By inductional hypothesis we have
\beaa
I(\pi,\t)\!\!\!&=&\!\!\! \mathbb{E}\Big\{\int_{t}^{t_1}U(X^{t,x}_{s},Q^{\pi,q}_{s})d\pi_{s}+V(t_1,X^{t,x}_{t_1},\pi_{t_1},Q^{\pi,q}_{t_1})
{\bf 1}_{\{\t=t_1\}} \\
\!\!\!&&\!\!\!+ \hE\Big[\int_{t_1}^{\tau}U(X^{t,x}_{s},Q^{\pi,q}_{s})d\pi_{s}+V(\tau,X^{t,x}_{\tau},\pi_{\tau},Q^{\pi,q}_{\tau})\Big|\cF_{t_1}\Big]
{\bf 1}_{\{\t > t_1\}} \Big\} \\
\!\!\!&\ge&\!\!\!  \mathbb{E}\Big\{\int_{t}^{t_1}U(X^{t,x}_{s},Q^{\pi,q}_{s})d\pi_{s}+V(t_1,X^{t,x}_{t_1},\pi_{t_1},Q^{\pi,q}_{t_1}){\bf 1}_{\{\t=t_1\}} + V(t_1,X^{t,x}_{t_1},\pi_{t_1},Q^{\pi,q}_{t_1}){\bf 1}_{\{\t > t_1\}} \Big\} \\
\!\!\!&=&\!\!\! \mathbb{E}\Big\{\int_{t}^{t_1}U(X^{t,x}_{s},Q^{\pi,q}_{s})d\pi_{s}+V(t_1,X^{t,x}_{t_1},\pi_{t_1},Q^{\pi,q}_{t_1})  \Big\}  \ge V(t,x,k,q),
\eeaa
where the last inequality is due to Proposition \ref{prop-DPPProp}. Since $\pi \in\sA^c_{ad}(t,k,q)$ is arbitrary, we proved \reff{DPPineq}  for $m$, completing the induction.

To prove \reff{DPPineq} for arbitrary $\t \in \cT_t$, we first find $\t_n \in \cT_t$, $n=1,2,\cds$, such that $\t_n - \t \le {1\over n}$ and 
$\t_n \downarrow \t$, as $n\to\infty$. By previous arguments we see that \reff{DPPineq} holds for each $\t_n$.  That is, $V(t,x,k,q) \le I(\pi, \t_n)$  for each $\pi \in\sA^c_{ad}(t,k,q)$. Moreover, by definition of $I(\pi, \t)$ we have 
\beaa
I(\pi, \t_n) - I(\pi,\t) \!=\! \hE\left\{\int_{\t}^{\tau_n}U(X^{t,x}_{s},Q^{\pi,q}_{s})d\pi_{s}+V(\tau_n,X^{t,x}_{\tau_n},\pi_{\tau_n},Q^{\pi,q}_{\tau_n}) -V(\tau,X^{t,x}_{\tau},\pi_{\tau},Q^{\pi,q}_{\tau}) \right\}.
\eeaa
Applying Corollary \ref{cor-Vregt} and noting that $\pi$ is continuous we see that the right hand side above converges to $0$ as 
$n\to\infty$. Consequently we obtain that $V(t,x,k,q) \le I(\pi, \t)$  for each $\pi \in\sA^c_{ad}(t,k,q)$. This implies \reff{DPPineq}, 
and hence concludes the proof. 
\qed

\begin{rem}
\label{rem-DPP}
{\rm Combining Theorems \ref{thm-DPP} and \ref{thm-equiV0V1}, we have the following alternative version of dynamic programming principle corresponding to the cost functional $J^1$ defined in \reff{Vt}:
\bea
\label{DPP2}
V(t,x,k,q)=\dneg \inf_{\pi\in\sA_{ad}(t,k,q)}\neg
\mathbb{E}\Big\{\dneg \int_{t}^{\tau}\neg U(X_{s},Q^{\pi}_{s})d\pi^c_{s}+ \dneg\sum_{t\le s<\t}\dneg D(X_s, Q^\pi_s, \D \pi_s)\neg+\neg  V(\tau,X_{\tau},\pi_{\tau},Q^{\pi,q}_{\tau})\Big\}. 
\eea
}
\end{rem}

\section{The HJB equation}
\setcounter{equation}{0}
\label{sect-HJB}

In this section we shall prove that
the value function, while not necessarily smooth, is a viscosity solution of the Hamilton-Jacobi-Bellman equation of
 the optimal execution problem. 

We begin by introducing some notations. For simplicity we often use the equivalent notations for partial derivatives: $\pa_t\vf=\frac{\pa\vf}{\pa t}$. The notations  $\pa_x\vf$, $\pa_k\vf$, $\pa_q\vf$, and $\pa_{xx}\vf$ are thus obvious. In this and next section, we denote by $C^{1,2}_b([0, T]\times \bar\cO)$ the set of continuous functions $\vf$ on $[0, T]\times \bar\cO$ such that the partial derivatives $\pa_t
\vf$, $\pa_x\vf$, $\pa_k\vf$, $\pa_q\vf$, and $\pa_{xx}\vf$ exist and are continuous and bounded. For each $t\in [0, T)$, we introduce a new filtration:
\bea
\label{hatdbF}
\hat \hF^t := \{\hat \cF^t_s\}_{s\ge 0} := \{\cF^W_s \vee \cF^Y_{s\wedge t}\}_{s\ge 0}.
\eea
Moreover,  in light of the cost functional $J^1$  in \reff{Vt} and the DPP (\ref{DPP2}), 
we define, for each $(t,x,k,q) \in [0, T)\times \bar\cO$, $\pi \in \sA_{ad}(t,x,k,q)$, $\vf\in C([0, T]\times \bar\cO)$, and $\hF$-stopping time $\t$,  
\bea
\label{Ipi}
I(\vf, \pi, \t) \!:=\!\mathbb{E}\Big\{\int_{t}^{\t} U(X_{s},Q^{\pi}_{s})d\pi^c_{s}+\neg \sum_{t\le s< \t}\dneg D(X_s, Q^\pi_s, \D \pi_s)+  \vf(\tau,X_{\tau},\pi_{\tau},Q^{\pi}_{\tau})\Big\} 
- \vf(t,x,k,q).
\eea

Next, we let $\t^t_1$ be  the first jump time of $N$ after $t$ and $\nu$ is   the common distribution of the jump size
random variables $\L_i$'s. We remark here that, by definition \reff{hatdbF} it is clear that $(\t^t_1, \D Y_{\t^t_1})$ is independent of $\hat \hF^t$, and hence $\t^t_1$ is not an $\hF^t$-stopping time(!). Furthermore, we have the following result that is important for our 
discussion.

\begin{lem}
\label{lem-0}
For any fixed $(t,k, q)$ and any $\pi\in\sA_{ad}(t,k,q)$, there exists an $\hat\hF^t$-adapted process $\tilde \pi$ such that $\tilde \pi_{s\wedge \t^t_1} = \pi_{s\wedge  \t^t_1}$, for all $s\ge t$, $\hP$-a.s.
\end{lem}

{\it Proof.}  We first note that since $\pi$ is left continuous, we need only find a $\hat\hF^t$-adapted process $\tilde\pi$ such that,
for any fixed $s\ge t$   $\hP\{\tilde \pi_{s}{\bf 1}_{\{\t^t_1>s\}} = \pi_{s}{\bf 1}_{\{\t^t_1>s\}}\}=1$. This amounts to saying that
given $s\ge t$, and $X\in \hL^0(\cF_s)$, there exists $\tilde X\in \hL^0(\hat\cF^t_s)$
such that 
$X{\bf 1}_{\{\t^t_1>s\}}=\tilde X{\bf 1}_{\{\t^t_1>s\}}$, $\dbP$-a.s. But this last statement is more or less standard (see, e.g., \cite{Bieleckis2004}), 
 we nevertheless give a brief proof for completeness. We fix $s >t$ and denote
 $$\cH_s :=\{X\in \hL^0(\cF_s)~|~\exists \tilde X\in \hL^0(\hat\cF^t_s), \mbox{ such that } X{\bf 1}_{\{\t^t_1>s\}}=\tilde X{\bf 1}_{\{\t^t_1>s\}}, ~\dbP\mbox{-a.s.}\}.
 $$
 Clearly, $\cH_s\subseteq \hL^0(\cF_s)$. We claim that $\cH_s\supseteq  \hL^0(\cF_s)$. Indeed, note that $\cF_s=\hat\cF^t_s\vee \si\{Y_r, t\le r\le s\}$. By a simple Monotone Class argument, for any $X\in \hL^0(\cF_s)$, we need only assume either $X\in \hL^0(\hat\cF^t_s)$ or $X=Y_r$ 
 for some $r\in [t,s]$. But in the former case we can
 choose $\tilde X= X$, and in the latter case we choose $\tilde X=Y_t$. Since in both cases $\tilde X\in \hL^0(\hat\cF^t_s)$, we conclude that $X\in \cH_s$. This 
 proves the claim, whence the lemma.
 \qed

Now for any $\vf \in C^{1,2}_b([0, T]\times \bar\cO$ we introduce the following integro-differential operators:
\bea
\label{cLcM}
 \sL[\vf](t,x,k,q) &:=& (\pa_t\vf+b\pa_x\vf+\frac{1}{2}\sigma^{2}\pa_{xx}\vf)(t,x,k,q) \nonumber\\
 &&+ \l  \int_{\hR}\big[\vf(t,x,k, (q+u)^+)-\vf(t,x,k,q)\big]\n(du);\\
\sM[\vf](t,x,k,q) &:=& U(x,q)+(\pa_k\vf-\pa_{q}\vf)(t,x,k,q).\nonumber
\eea
The following lemma is crucial.
\begin{lem}
\label{lem-I}
Assume $\vf \in C^{1,2}_b([0, T]\times \bar\cO$ and $\t$ is an $\hF^t$-stopping time. Then it holds that
\bea
\label{Ipi1}
I(\vf, \pi, \t \wedge \t^t_1) &=& \mathbb{E}\Big\{\int_t^{\t \wedge \t^t_1}  \sL[\vf](s,  X_s, \pi_s,  Q^\pi_s) ds + \int_t^{\t \wedge \t^t_1} \sM[\vf] (s,  X_s, \pi_s,  Q^\pi_s) d\pi^c_s \nonumber\\
&&+\sum_{t\le s<\t \wedge \t^t_1}  \int_0^{\D \pi_s} \sM[\vf](s,  X_s, \pi_s + u, Q^\pi_s - u) du\Big\}.
\eea
where $\sL$ and $\sM$ are  defined by (\ref{cLcM}).
\end{lem}

{\it Proof.}  
For any $\hF^t$-stopping time $\t$ we denote
$\hat\t := \t \wedge \t^t_1$. Let $\pi\in\sA_{ad}(t,k,q)$, and let $\tilde \pi$ be the $\hF^t$-adapted version of $\pi$ defined
in Lemma \ref{lem-0}, and define
$\tilde Q^\pi_s :=  q-\tilde \pi_s + k$, $s\ge t$.
Then, it is readily seen that $Q^{\pi}_{\tau^t_1}= (\tilde Q^{\pi}_{\tau^t_1} + \D Y_{\t^t_1})^+$, and thus
\bea
\label{Iest1}
&&\vf(\hat\tau,X_{\hat\tau},\pi_{\hat\tau},Q^{\pi}_{\hat\tau})  - \vf(t,x,k,q) =\vf(\hat\tau,X_{\hat\tau},  \pi_{\hat\tau},\tilde Q^{\pi}_{\hat\tau})  - \vf(t,x,k,q)\\
&&\qq \qq \qq\qq+  \Big[\vf(\tau^t_1,X_{\t^t_1}, \tilde\pi_{\t^t_1}, (\tilde Q^{\pi}_{\t^t_1} + \D Y_{\t^t_1})^+) - \vf(\tau^t_1,X_{\t^t_1}, \tilde\pi_{\t^t_1}, \tilde Q^{\pi}_{\t^t_1})\Big]{\bf 1}_{\{\t^t_1 \le \t\}}.\nonumber
\eea
 Since $(X, \tilde \pi, \tilde Q^\pi)$, $\t^t_1$, $\D Y_{\t^t_1}$ are independent, we have
 \bea
 \label{Iest2}
&& \hE\Big\{ \Big[\vf(\tau^t_1,X_{\t^t_1}, \tilde\pi_{\t^t_1}, (\tilde Q^{\pi}_{\t^t_1} + \D Y_{\t^t_1})^+) - \vf(\tau^t_1,X_{\t^t_1}, \tilde\pi_{\t^t_1}, \tilde Q^{\pi}_{\t^t_1})\Big]1_{\{\t^t_1 \le \t\}}\Big\}\nonumber\\
&=&\hE\Big\{ \int_t^\t \l e^{-\l(s-t)} ds \int_{\hR} \Big[\vf(s,X_s, \tilde\pi_s, (\tilde Q^{\pi}_s + u)^+) - \vf(s,X_s, \tilde\pi_s, \tilde Q^{\pi}_s)\Big]\n(du)\Big\}\\
&=& \hE\Big\{\l \int_t^{\t} 1_{\{\t^t_1\ge s\}} \int_{\hR} \Big[\vf(s,X_s, \tilde\pi_s, (\tilde Q^{\pi}_s + u)^+) - \vf(s,X_s, \tilde\pi_s, \tilde Q^{\pi}_s)\Big]\n(du)ds\Big\}\nonumber\\
&=& \hE\Big\{\l \int_t^{\hat\t}  \int_{\hR} \Big[\vf(s,X_s, \pi_s, ( Q^{\pi}_s + u)^+) - \vf(s,X_s, \pi_s,  Q^{\pi}_s)\Big]\n(du)ds\Big\}.\nonumber
  \eea
Here  we used the fact that $\tilde Q^\pi_s = Q^\pi_s$, $t\le s<\hat\t$. Furthermore, applying It\^{o}'s formula we have
\bea
\label{Iest3}
&&\hE\Big\{\vf(\hat\tau,X_{\hat\tau},  \pi_{\hat\tau},\tilde Q^{\pi}_{\hat\tau})  - \vf(t,x,k,q)\Big\}\nonumber\\
&=&\hE\Big\{\int_t^{\hat\t} \Big[\pa_t\vf+b\pa_x\vf+\frac{1}{2}\sigma^{2}\pa_{xx}\vf\Big](s, X_s, \pi_s,  Q^\pi_s) ds + \int_t^{\hat\t} \Big[\pa_k\vf - \pa_q\vf\Big](s, X_s, \pi_s,  Q^\pi_s) d\pi^c_s\nonumber\\
&&+\sum_{t\le s<\hat\t}  \int_0^{\D \pi_s} \Big[\pa_k \vf- \pa_q\vf\Big](s,  X_s, \pi_s + u,  Q^\pi_s - u) du  \Big\}.
\eea
Plugging \reff{Iest2}, \reff{Iest3} into \reff{Iest1}, and then plugging \reff{Iest1}, \reff{CgeD} into \reff{Ipi}, we obtain \reff{Ipi1}.
\qed

It is worth noting that if we use the continuous strategy $\pi \in \sA^c_{ad}(t,x,k,q)$, then \reff{Ipi} and \reff{Ipi1} become
\bea
\label{Ipi2}
I(\vf, \pi, \t) &:=& \mathbb{E}\Big\{\int_{t}^{\t} U(X_{s},Q^{\pi}_{s})d\pi_{s}+  \vf(\tau,X_{\tau},\pi_{\tau},Q^{\pi}_{\tau})\Big\} 
- \vf(t,x,k,q)\\
\label{Ipi3}
&=& \mathbb{E}\Big\{\int_t^{\t}  \sL[\vf] (s,  X_s, \pi_s, Q^\pi_s) ds + \int_t^{\t} \sM[\vf] (s,  X_s, \pi_s, \tilde Q^\pi_s) d\pi_s\Big\},
\eea
respectively. Clearly,  \reff{Ipi2} is valid even when $\vf$ is not smooth. In fact,  by Theorem \ref{thm-DPP} we have
\bea
\label{IV}
0 &=& \inf_{\pi\in\sA^c_{ad}(t,x,k,q)} I(V, \pi, \t).
\eea

Furthermore,  if $V\in C^{1,2}_b([0, T]\times \bar \cO)$, then we may plug \reff{Ipi3} into \reff{IV} and deduce  
 the following Quasi-Variational-Inequality (QVI):
\bea
\label{HJB}
    \min\Big(\sL[V], ~\sM[V]\Big)(t,x,k,q)=0,\q (t,x,k,q)\in [0, T)\times \cO,
 \eea
with the terminal-boundary conditions:
\bea
\label{terminal}
        V(T,x,k,q)=g(x,K-k);\q
    V(t,x,K,q)=0;\q
\sL[V](t,x,k,0)=0.
\eea
As we will see in next section, in this case $V$ is indeed the unique classical solution of the QVI (\ref{HJB}) and (\ref{terminal}).

In general, however, $V$ may not be smooth. We thus need to make use of  the notion of the {\it viscosity solution}. To this end, 
let us denote, 
for $(t,x,k,q)\in [0, T) \times \hR_+ \times [0, K) \times \bar\hR_+$,
\bea
\label{cA}
\sA(t,x,k,q) \neg&\neg:=\neg&\neg \Big\{\vf\in C^{1,2}_b([0,T]\times \bar\cO):  [V-\vf](t,x,k,q)=0\Big\};\nonumber\\
\overline\sA(t,x,k,q) \neg&\neg:=\neg& \neg\Big\{\vf\in \cA(t,x,k,q):   \mbox{ $V-\vf$ attains a global maximum at
$(t,x,k,q)$}\Big\} ;\\
 \underline\sA(t,x,k,q) \neg&\neg:=\neg& \neg\Big\{\vf\in \cA(t,x,k,q):   \mbox{ $V-\vf$ attains a global minimum at
$(t,x,k,q)$}\Big\}.\nonumber 
\eea

\begin{defn}
\label{defn-visdef}
A continuous function $V: [0, T]\times \bar\cO\mapsto \mathbb{R}_+$
is called a viscosity subsolution (resp. supersolution) to the QVI  (\ref{HJB})-(\ref{terminal}) if 

(i) $V(T,x,k,q)\ge (resp. \le ) g(x,K-k)$ and $V(t,x,K,q) \ge 0 (resp. \le 0)$;

(ii) for any $(t,x,k,q)\in [0, T)\times \cO$ and $\vf\in \overline\cA(t,x,k,q)$ (resp. $\underline\cA(t,x,k,q)$)  one has:
$$min(\sL[\vf],\sM[\vf])(t,x,k,q)\geq0, (resp. \leq0);$$

(iii) for any $(t,x,k)\in [0, T)\times \hR_+ \times [0, K)$ and $\vf\in \overline\cA(t,x,k,0)$ (resp. $\underline\cA(t,x,k,0)$)  one has:
$$\sL[\vf](t,x,k,0)\geq0, (resp. \leq0).$$

\no Moreover, $V$ is called a viscosity solution if it is both a viscosity subsolution and  supersolution.
\qed
\end{defn}

Our main result of this section is the following theorem.
\begin{thm}
\label{thm-viscosity}
Assume (H1)-(H3). Then the value
function $V$ of the optimal execution problem is a viscosity
solution of the QVI \reff{HJB}-\reff{terminal}.
\end{thm}

{\it Proof.} The terminal condition $V(T,x,k,q)=g(x,K-k)$ is obvious. Moreover, note that if $\pi_t=K$, then $\pi_s\equiv K$ for all $s\in[t,T)$, as there is no need to purchase any more. Thus
$d\pi_s=0$ for $s\in[t,T]$, and clearly $g(X_T,K-\pi_T)=g(X_T,0)=0$. That is,  $V(t,x,K,q)=0$. So Definition \ref{defn-visdef} (i)  holds (with equalities), and thus it suffices to check Definition \ref{defn-visdef} (ii) and (iii). 

We first prove the viscosity subsolution properties.  It suffices to show that, for any  $(t,x,k,q)\in [0, T) \times \hR_+ \times [0, K) \times \bar\hR_+$ and $\vf\in \overline\cA(t,x,k,q)$. 
\bea
\label{subvis}
\sL[\varphi](t,x,k,q)\geq0,~\mbox{for}~q\ge 0; \qq
\sM[\varphi](t,x,k,q)\geq0,~\mbox{for}~q>0.
\eea

In what follows we denote, for 
$\d>0$ small,  $\t_\d :=(t+\d) \wedge \t^t_1$, and let $C>0$ be a generic constant that is allowed to vary from line to line.

We begin by proving the first inequality in \reff{subvis}. Let $\pi := k$ be the constant process. Then $Q^\pi_s = q$ for $t\le s<\t_\d$.  By \reff{IV}, \reff{cA}, and \reff{Ipi3}, we have 
\bea
\label{DPPEq1}
0 &\le& I(V,k,\t_\d) \le I(\vf, k,\t_\d) = \mathbb{E}\Big\{\int_t^{\t_\d}  \sL[\vf] (s,  X_s, k, q) ds  \Big\} \nonumber\\
&=& \mathbb{E}\Big\{\int_t^{t+\d}   \sL[\vf] (s,  X_s, k, q) ds  \Big\} -  \mathbb{E}\Big\{\int_{\t^t_1}^{t+\d}   \sL[\vf] (s,  X_s, k, q) ds  1_{\{\t^t_1 < t+\d\}}\Big\}.
\eea
Note that $\sL[\vf] $ is bounded and 
\bea
\label{tauest1}
\hP(\t^t_1 < t+\d) \le C\d, 
\eea
dividing both sides of \reff{DPPEq1} by $\d$ and then sending $\d\to 0$, we prove the first inequality in \reff{subvis}.

To check the second inequality  in \reff{subvis} for $q>0$, let $\eta>0$ and set $\pi_s := k + {(s-t)\wedge \d \over \d} \eta q$. Clearly $\pi\in\sA^c_{ad}(t,k,q)$, $\tilde \pi = \pi$, $d \pi_s = {\eta q \over \d} ds$, and $\tilde Q^\pi_s = [1-{s-t\over \d} \eta] q$, $s\le \t_\d$.  By \reff{IV}, \reff{cA}, and \reff{Ipi3} again, we have  
\beaa
0&\le&\mathbb{E}\Big\{{\eta q\over \d} \int_{t}^{\t_\d}\sM
[\vf](s,X_s,\pi_{s},[1-{s-t\over \d} \eta] q)d{s}+\int_{t}^{\t_\d}\sL
[\vf](s,X_s,\pi_{s},\tilde Q^\pi_s)ds\Big\}\\
&\le&\mathbb{E}\Big\{{\eta q\over \d} \int_{t}^{t+\d}\sM
[\vf](s,X_s,\pi_{s},[1-{s-t\over \d} \eta] q)d{s}\Big\} + C \hP(\t^t_1 < t+\d) +C\d \\
&\le& \mathbb{E}\Big\{{\eta q\over \d} \int_t^{t+\d} \big[\sup_{0\le \th \le 1}  \sM
[\vf](s,X_s,\pi_{s},[1-\th \eta] q)\big] ds\Big\} + C\d,
\eeaa
Here in the last inequality above we used \reff{tauest1} again. Now, sending $\d\to 0$ in the above
we can easily deduce that $ \sup_{0\le \th \le 1}\sM [\vf](t, x, k,[1-\th\eta ] q)\ge 0$. The arbitrariness of $\eta>0$ then further leads
to the second inequality of \reff{subvis},   proving the viscosity subsolution property.

 We now turn to the viscosity supersolution property.  We first check Definition \ref{defn-visdef} (iii). Let  $(t,x,k)\in [0, T)\times \hR_+ \times [0, K)$ and  $\vf\in\underline\cA(t,x,k,0)$. For any $\pi \in \sA^c_{ad}(t,k,0)$, since there is no
liquidity ($q=0$), there is no possibility of trading, and thus it must hold that: $\pi_s\equiv k$ and $Q^{\pi,0}_s = 0$, $s<\t^t_1$. Then, by   \reff{IV}, \reff{cA} and \reff{Ipi3} again, we have
\beaa
0 = I(V, k, \t_\d) \ge I(\vf,  k, \t_\d) = \mathbb{E}\Big\{\int_{t}^{\t_\d} \sL [\vf](s,X_{s},k, 0)ds \Big\}.
\eeaa
Dividing both sides above by $\d$ and then sending $\d\to 0$, similar to the case \reff{DPPEq1} we can prove Definition \ref{defn-visdef} (iii). 

It remains to verify Definition \ref{defn-visdef} (ii). Suppose in the contrary that 
\bea
\label{c}
c &:=& \min\Big(\sL[\vf],\sM[\vf]\Big)(t,x,k,q) > 0
\eea
 for some $(t,x,k,q)\in [0, T)\times \cO$ and $\vf\in \underline\cA(t,x,k,q)$.  Then, applying Theorem \ref{thm-DPP} on $\t^t_1$ we 
 can find $\pi:= \pi^\d\in \sA^c_{ad}(t,k,q)$ such that
\beaa
V(t,x,k,q) \ge  \hE\Big\{\int_{t}^{\t^t_1}U(X^{t,x}_{s},Q^{\pi}_{s})d\pi_{s}+V(\t^t_1, X^{t,x}_{\t^t_1}, \pi_{\t^t_1}, Q^{\pi, q}_{\t^t_1})\Big\} -\d^2.
\eeaa
Now let $\tilde \pi$ be the $\hF^t$-adapted version of $\pi$, as was defined in Lemma \ref{lem-0}, and $\tilde Q^\pi_s=q-\tilde \pi_s+k$, $s\ge t$. For any $\d>0$, define the following stopping times:
\bea
\label{td}
\t^X_\d &:=&  \inf\big\{s>t: |X^{t,x}_s - x|\ge \d^{1\over 4}\big\}\wedge T,\q \t^\pi_\d := \inf\big\{s>t: \tilde\pi_s - k \ge \d\big\}\wedge T, \nonumber \\
\t'_\d &:=&  (t+\d)\wedge \t^X_\d \wedge \t^\pi_\d,\q \hat \t'_\d :=\t'_\d \wedge \t^t_1.
\eea
Then $\t'_\d$ is an $\hat\hF^t$-stopping time.  Similar to the first part of Proposition \ref{prop-DPPProp} we can show
that
\beaa
 V(t,x,k,q) &\ge&  \hE\Big\{\int_{t}^{\hat \t'_\d}U(X_{s},Q^{\pi}_{s})d\pi_{s}+V(\hat \t'_\d, X_{\hat \t'_\d},  \pi_{\hat \t'_\d}, Q^{\pi}_{\hat \t'_\d})\Big\} -\d^2.
\eeaa
Now following the derivation of  \reff{DPPEq1} we obtain
\bea
\label{d2}
\d^2 &\ge& \mathbb{E}\Big\{\int_{t}^{\hat \t'_\d}\sM
[\vf](s,X_s, \pi_{s}, Q^{\pi}_s)d \pi_s+\int_{t}^{\hat\t'_\d} \sL [\vf]  (s,X_s, \pi_{s},Q^{\pi}_s)ds\Big\}.
\eea
Since $\vf$ is smooth,  we deduce from \reff{c}  that, for $\d$ is small enough, 
 \beaa
\sM[\vf](s,X_s, \pi_{s}, Q^{\pi}_s)  \ge {c\over 2},\q (\sL [\vf] + \sG[\vf]\big)(s,X_s,\pi_{s}, Q^{\pi}_s) \ge {c\over 2},\q t\le s <\hat \t'_\d.
\eeaa
Thus it follows from (\ref{d2}) that 
$ \dis \d^2 \ge {c\over 2} \mathbb{E}\{  \pi_{\hat\t'_\d} - k + \hat\t'_\d - t  \}$.
But note that $\pi_{\hat\t'_\d} - k = \d$ on $\{\hat\t'_\d = \t^\pi_\d\}$, this leads further to
\bea
\label{d2a}
\d^2 &\ge& {c\over 2} \mathbb{E}\big\{  \d {\bf 1}_{\{\hat\t'_\d = \t^\pi_\d\}} + ((t+\d)\wedge \t^X_\d \wedge \t^t_1 - t)  1_{\{\hat\t'_\d < \t^\pi_\d\}} \big\}\ge {c \over 2} \d - C \mathbb{E}\big\{  (t+\d - \t^X_\d \wedge \t^t_1)^+ \big\}\nonumber\\
&\ge& {c  \over 2} \d - C \d \big[ \hP(\t^X_\d < t+\d) + \hP(\t^t_1 < t+\d)\big].
\eea
Finally, recalling \reff{tauest1} and noting that
\beaa
\hP(\t^X_\d  < t+\d)  = \hP\Big(\sup_{t\le s\le t+\d} |X^{t,x}_s - x| \ge \d^{1\over 4}\Big) \le {1\over \d}  \hE\Big\{\sup_{t\le s\le t+\d} |X^{t,x}_s - x|^4 \Big\} \le C(1+|x|^4) \d.
 \eeaa
 We derive from (\ref{d2a}) that
 $
 \d^2  \ge {c\over 2}\d - C(1+|x|^4)\d^2
 $. But this is obviously impossible  when $\d>0$ is small enough, a contradiction to the assumption (\ref{c}). This completes the proof.
 \qed

 \begin{rem}
 \label{rem-unique}
 {\rm (i) If the value function actually has the regularity $V\in C^{1,2}_b([0, T]\times \bar \cO)$, then instead of being
 a viscosity solution, it will be a classical solution to the QVI (\ref{HJB}).  Moreover, by Theorem \ref{thm-existence} below, we see that the classical solution is unique.
 
 (ii) We should note that one may try to analyze the uniqueness in the sense of  viscosity solutions by following the
 standard techniques (see the classical reference \cite{CIL}). However, since our main focus is the  dynamic equilibrium model of the limit order book, we prefer not to pursue this 
 in this already lengthy paper and will leave it to interested reader.
 \qed}
 \end{rem}
 
\section{Description of Optimal Strategy}
\setcounter{equation}{0}
\label{sect-strategy}

In this section we give a characterization of the optimal strategy. Our argument will be based on the assumption that
the HJB equation has a ``classical solution", which will not be substantiated in this paper, as it is itself a challenging problem.
Our main purpose is to see the possible structure of the optimal strategy and compare it to the usual optimal singular stochastic
control in the literature. 

Our starting point is the following partial {\it Verification Theorem}.
\begin{prop}
\label{prop-verification}
Assume (H1) - (H3), and that $v\in C^{1,2}_b([0, T]\times \bar \cO)$ is a classical solution to the QVI (\ref{HJB})-\reff{terminal}.  Then $v\le V$.
\end{prop}

{\it Proof.} Without loss of generality, we assume $t=0, k=0$.  By \reff{Vt}, it suffices to show that
\bea
\label{verif}
v(0,x,0,q)\leq J^1(0,x,0,q;\pi), &\mbox{for any}& \pi\in\sA_{ad}(0,0,q).
\eea
We remark that, for this proposition, we can actually utilize $J^0$, namely considering only continuous strategies. However, to analyze the optimal strategy later, we shall use $J^1$ instead.

Recall that $0<\t_1<\t_2,\cds$ are the jump times of $N$.  Denote $\hat\t_i := \t_i \wedge T$. By the terminal condition  \reff{terminal},  we have
\beaa
e_\pi &:=& J^1(0,x,0,q;\pi) -  v(0,x,0,q)  \\
&=&\mathbb{E}\Big\{\int_0^{T}U(X_{s},Q^{\pi}_{s})d\pi^c_s+\sum_{0\leq s<T}D(X_{s},Q^\pi_{s},\Delta\pi_{s})+ v(T, X_T, \pi_T, Q^\pi_T) - v(0, X_0, \pi_0, Q^\pi_0)\Big\}\\
&=&  \sum_{i=0}^\infty \mathbb{E}\Big\{\int_{\hat\t_i}^{\hat\t_{i+1}}U(X_{s},Q^{\pi}_{s})d\pi^c_s+\sum_{\hat\t_i\leq s<\hat\t_{i+1}}D(X_{s},Q^\pi_{s},\Delta\pi_{s})\\
&&\qq + v(\hat\t_{i+1}, X_{\hat\t_{i+1}}, \pi_{\hat\t_{i+1}}, Q^\pi_{\hat\t_{i+1}}) - v(\hat\t_i, X_{\hat\t_i}, \pi_{\hat\t_i}, Q^\pi_{\hat\t_i})\Big\}.
\eeaa
By introducing the filtrations  $\hat\dbF^i := (\cF^W_s \vee \cF^Y_{s\wedge \hat\t_i})_{0\le s\le T}$ and setting $\t:= T$ in \reff{Ipi1}, we obtain
\bea
\label{epi}
e_\pi& = &\sum_{i=0}^\infty  \mathbb{E}\Big\{\int_{\hat\t_i}^{\hat\t_{i+1}}  \sL[v](s,  X_s, \pi_s,   Q^{\pi}_s) ds + \int_{\hat\t_i}^{\hat\t_{i+1}} \sM[v] (s,  X_s,  \pi_s,  Q^{\pi}_s) d\pi^c_s \nonumber\\
&& +\sum_{\hat\t_i\le s<\hat\t_{i+1}}  \int_0^{\D \pi_s} \sM[v](s,  X_s,  \pi_s + u, Q^{\pi}_s - u) du \Big\} \ge 0,
\eea
thanks to \reff{HJB}. This completes \reff{verif}.
\qed

In the rest of the section we shall find an optimal strategy $\pi^*\in\sA_{ad}(0,0,q)$ such that \reff{epi}, hence \reff{verif}, holds with equality, given the existence of the classical solution $v$ of the QVI \reff{HJB}-\reff{terminal}. We shall remark though, although it is interesting in theory, the $\pi^*$ is in general not implementable since the cost $D$ in the expression $J^1$ of \reff{Vt} is not the real jump cost. However, as was pointed out in Remark \ref{rem-J01}, this $\pi^*$ will nevertheless provide us  a very good and implementable approximate optimal strategy.

To help identifying the optimal strategy $\pi^*$, we first provide some sufficient conditions. Without loss of generality, we shall only focus on the interval $[0, \t_1]$, corresponding to the term in \reff{epi} with $i=0$. To be more precise, we want to find $\pi\in\sA_{ad}(0,0,q)$ such that
\bea
\label{epi0}
e_{\pi,0}& := &\mathbb{E}\Big\{\int_{0}^{\hat\t_1}  \sL[v](s,  X_s, \pi_s,  Q^{\pi}_s) ds + \int_{0}^{\hat\t_1} \sM[v] (s,  X_s,  \pi_s,  Q^{\pi}_s) d\pi^c_s \nonumber\\
&& +\sum_{0\le s<\hat\t_{1}}  \int_0^{\D \pi_s} \sM[v](s,  X_s,  \pi_s + u, Q^{\pi}_s - u) du \Big\} = 0.
\eea
To this end, for any $(t,x,k,q)\in[0,T]\times\hR_+\times [0, K] \times \bar\hR_+$, denote
\bea
\label{Of}
\left.\ba{lll}
\dis O(t,x,q) &:=& \big\{ y\in [0, K\wedge q]:\sM[v](t,x, y,q-y)>0\big\} ;\\
\dis \f(t,k,q) &:=& \inf\big\{y>k: y\in O(t, X_t, q)\big\} \wedge K \wedge q.
\ea\right.
\eea
It is clear that $O(t,x,q)$ is an open set in $[0, K\wedge q]$, and $\f$ is  $\hF^W$- progressively measurable, 
non-decreasing in $k$, such that $\f(t,k, q)\geq k$,  and  $\f(t,k,q)=k$ for  $k\in O(t,X_t,q)$.
We  have the following result.
\begin{prop}
\label{prop-suff}
Assume all the conditions of Proposition \ref{prop-verification} hold. If $\pi\in\sA_{ad}(0,0,q)$ satisfies:
\bea
\label{suff}
\int_{0}^{\hat\t_{1}}\mathbf{1}_{O(t,X_{t}, q)}(\pi_{t})d\pi^c_{t}=0 &\mbox{and}&
\pi_{t+}=\phi(t,\pi_{t}, q), ~t\in [0,\hat\t_1), \q\hP\mbox{-a.s.}
\eea
then \reff{epi0} holds.
\end{prop}

{\it Proof.} First,  denote $O^c (t,x,q) := [0, K\wedge q] - O(t,x,q)$. Then the first equality in \reff{suff} implies: 
\beaa
d\pi^c_t= \Big[\mathbf{1}_{O(t,X_{t}, q)}(\pi_{t}) + \mathbf{1}_{O^c(t,X_{t}, q)}(\pi_{t})\Big]d\pi^c_{t}  = \mathbf{1}_{O^c(t,X_{t}, q)}(\pi_{t}) d\pi^c_t,\q 0\le t\le \hat\t_1.
\eeaa
Note that $Q^\pi_t = q-\pi_t$, $0\le t<\hat\t_1$, then by the definition of $O$ in  \reff{Of} we have
\bea
\label{epi1}
 \int_{0}^{\hat\t_1} \sM[v] (s,  X_s,  \pi_s,  Q^{\pi}_s) d\pi^c_s=   \int_{0}^{\hat\t_1} \sM[v] (s,  X_s,  \pi_s,  Q^{\pi}_s)  \mathbf{1}_{O^c(t,X_{s}, q)}(\pi_{s})d\pi^c_s=0.
 \eea

Next,  when $\D \pi_t >0$, by the second condition of \reff{suff}  we have
\beaa
\pi_{t+}= \phi(t,\pi_{t}, q)=\inf\big\{y>\pi_{t}:\sM[v](t,X_{t}, y, q-y)>0\big\}\wedge K\wedge q.
\eeaa
This implies that $\sM[v](t,X_{t}, y, q-y) = 0$ for all $\pi_t \le y < \pi_{t+}$. Thus, by denoting $y = \pi_t + u$,
\bea
\label{epi2}
 \int_0^{\D \pi_s} \sM[v](s,  X_s,  \pi_s + u, Q^{\pi}_s - u) du =  \int_0^{\D \pi_s} \sM[v](s,  X_s,  \pi_s + u, q-\pi_s - u) du=0.
 \eea

Finally, we claim that
\bea
\label{LG=0}
\sL[v](t,X_t,\pi_t,q-\pi_t)=0 &\mbox{for}& t\in [0, \hat\t_1]~\mbox{such that}~\D \pi_t=0.
\eea
We note that if (\ref{LG=0}) is substantiated, then since $\pi$ has at most countably many jumps, we have
 \bea
 \label{epi3}
 \hE\Big\{\int_0^{\hat\t_1} \sL[v](t,X_t,\pi_t,Q^\pi_t) dt \Big\}= 0.
 \eea
 Combining \reff{epi1}, \reff{epi2}, and \reff{epi3}, we prove \reff{epi}.
 
 It remains to prove \reff{LG=0}.  Fix $ t\in [0, \hat\t_1]$ such that $\D \pi_t=0$. If $\pi_t=q$, then \reff{LG=0} is the third  condition of \reff{terminal}.  If $\pi_t = K$,  then $\pi_s = K$ for all $s\in [t,T]$, and thus $v(s, X_s, \tilde \pi_s, \tilde Q^\pi_s) = 0$, thanks to the second  condition of \reff{terminal}.  Compare \reff{Ipi2} and \reff{Ipi3}, one can easily check  \reff{LG=0}.  Now assume $\pi_t<K\wedge q$, then 
\beaa
\pi_{t}=\pi_{t+}=\phi(t,\pi_t, q)=\inf\big\{y>\pi_t: y\in O(t,X_{t}, q)\big\}.
\eeaa
That is,  $\pi_{t}\in\bar O(t,X_t, q)$.
But note that as the solution to the variational inequality (\ref{HJB}), it is easy to see that
$\sL[v](t,X_t,y,q-y)=0$  holds whenever $\sM[v](t, X_t,y, q-y)>0$, namely, for any $y\in O(t,X_t, q)$. 
The continuity of $v$ then renders  that $\cL[v](t,X_t,y,q-y)=0$ on 
$\bar{\cO}(t,X_{t}, q)$. Consequently, (\ref{LG=0}) holds. This proves \reff{epi3},
whence the theorem.
\qed

We next show that such $\pi$ indeed exists.  Fix $(x,q)$.  In light of Proposition \ref{prop-suff} we introduce:
\bea
\label{sA0}
\sA_{0}=\Big\{\pi\in\sA_{ad}(0,0,q): \int_{0}^{\hat\t_{1}}\mathbf{1}_{O(t,X_{t}, q)}(\pi_{t})d\pi^c_{t}=0,
\pi_{t+}\le \phi(t,\pi_{t}, q), ~t\in [0,\hat\t_1), ~\hP\mbox{-a.s.}
\Big\}.
\eea
Clearly,  $\pi_{t}\equiv0\in\sA_{0}$, thus $\sA_{0}\neq\emptyset$. We
shall construct the optimal strategy from this set. 
\begin{prop}
\label{prop-existence}
Assume all the conditions of Proposition \ref{prop-verification} hold. Then there exists $\pi\in\sA_0 \subset \sA_{ad}(0,0,q)$ satisfying \reff{suff}, and consequently  \reff{epi0} holds.
\end{prop}

{\it Proof.}  We shall prove the existence by using Zorn's lemma. To this end, we introduce a partial order in $\sA_{0}$:
\bea
\label{order}
\pi^{1}\prec \pi^{2} \q \mbox{  if and only if  }\q \pi^{1}_{t}\leq\pi^{2}_{t} \q\text{  for all  }t\in[0,T],\text{  $\hP$-a.s.}
\eea
We claim that every totally ordered subset in $\sA_0$ has an upper bound in $\sA_0$. Indeed, 
let  $ \{\pi^{i}\}_{i\in I}\subseteq\sA_{0}$ be a totally ordered subset, where the index set $I$ could be uncountable. 
Denoting $\hQ_T$ to be the set of all rationals in $[0,T]$,  we define
\bea
\label{suppi}
\pi_{r} := \esssup_{i\in I}\pi^{i}_{r}, \qq \forall r\in\mathbb{Q}_{T}.
\eea
Since $\{\pi^{i}\}$ is totally ordered, by a standard argument we can find a sequence $\pi^{n}=\pi^{i_n}$, $i_n\in I$, $n=1,2,\cds$, such
that $\pi^n$'s are 
non-decreasing in $n$; and 
\bea
\label{limpi}
\lim_{n}\pi^{n}_{r}=\esssup_{i\in I}\pi^{i}_{r}=\pi_{r}, \qq \forall r\in\mathbb{Q}_T.
\eea
We then define 
$\pi_{t}:=\lim_{r\nearrow t, r\in \hQ_T}\pi_{r}$, for all $t\in (0,T]$. We shall prove that $\pi\in\sA_{0}$, and therefore 
an upper bound of $\{\pi^{i}\}$. Clearly $\pi$ is $\hF$-adapted, non-decreasing, 
left continuous, and $\pi_0 = 0$, $\pi_T \le K$. Moreover, since $Q^{\pi^n}\ge 0$, clearly $Q^\pi_r\ge 0$ for all $r\in \mathbb{Q}_T$, which implies $Q^\pi_t \ge 0$ for all $t\in [0, T]$ and  thus $\pi\in\sA_{ad}(0,0,q)$. 

We now check that $\pi$ satisfies the two requirements of $\sA_0$. 
Since there is no stochastic integral involved, in what follows we shall  fix $\o\in\O$, modulo a $\hP$-null set, if necessary.

(i) We first show that $\int_0^{\hat\t_1}{\bf 1}_{O(t,X_t, q)}(\pi_{t})d\pi^c_t=0$.  Indeed, since $\pi$ has at most countably many jumps, it suffices to show that
\beaa
\int_0^{\hat\t_1}{\bf 1}_{O(t,X_t, q)}(\pi_{t}){\bf 1}_{\{\D \pi_t=0\}}d\pi^c_t=0.
\eeaa
Now for any $t\in [0, \hat\t_1)$ such that $\D \pi_t = 0$ and $\pi_t \in O(t,X_t,q)$, by \reff{Of} we have $\sM[v](t,X_{t},\pi_{t},q-\pi_{t})>0$. By the continuity of $\sM[v]$, there exists $\e>0$ such that 

(a) $\sM[v](s,X_s, y,q-y)>0$, for all $s\in [(t-\e) \vee 0, (t+\e)\wedge \hat\t_1]$; and 

(b) $y \in [(\pi_t -\e)\vee 0, (\pi_t+\e) \wedge K\wedge q]$. 

Since $\pi$ is continuous at $t$, there exists rationals $r_1, r_2$ such that $(t-\e) \vee 0 \le r_1 < t < r_2 \le (t+\e)\wedge \hat\t_1$ and $\pi_t - {\e\over 3} \le \pi_{r_1} \le \pi_t \le \pi_{r_2} \le \pi_t + {\e\over 3}$. 
Now by the monotone convergence of $\pi^n_r$, in the spirit of Dini's lemma, there exists $n_0$ such that, for all $n\ge n_0$, $|\pi^n_s - \pi_s|\le \e$ for $s\in [r_1, r_2]$. This implies $\sM[v](s,X_s, \pi^n_s,q-\pi^n_s)>0$, and thus $\pi^n_s \in O(s,X_s, q)$,  for all $s\in [r_1, r_2]$ and $n\ge n_0$. Since $\pi^n\in \sA_0$, then $\int_{r_1}^{r_2} d (\pi^n)^c_t = 0$ and $\pi^n_{s+} \le \f(s, X_s, \pi^n_s) = \pi^n_s$. That is, $\pi^n$ is a constant on $[r_1, r_2]$ for all $n\ge n_0$. Then $\pi$ is also a constant on $[r_1, r_2]$, and therefore, $\int_{r_1}^{r_2} {\bf 1}_{O(t,X_t, q)}(\pi_{t})d\pi^c_t=0$. Since $t$ is arbitrary, we prove the desired property.

(ii) We next show that $\pi_{t+}\leq\phi(t,\pi_{t}, q)$ for $t\in [0, \hat\t_1)$. For any $y\in (\pi_{t}, K\wedge q)$ such that $\sM[v](t,X_t, y,q-y)>0$. By the continuity of $\sM[v]$, there exists $0<\e < \hat\t_1-t$ such that  $\sM[v](s,X_s, y,q-y)>0$ for all $s\in [t, t+\e]$.  We claim that
\bea
\label{pinclaim}
\pi^n_s \le y,\q s\in [t, t+\e],~\mbox{for all}~n.
\eea
Note that if (\ref{pinclaim}) is true, then clearly $\pi_s \le y$ for $s\in [t, t+\e]$, which implies that $\pi_{t+} \le y$. By the arbitrariness of $y$, we obtain  $\pi_{t+}\leq\phi(t,\pi_{t}, q)$.

To see \reff{pinclaim}, suppose in the contrary that $\tilde t_n := \inf\{s\ge t: \pi^n_{s} > y\} <t+\e$. Then $\pi^n_{\tilde t_n} \le y \le  \pi^n_{\tilde t_n+}$. Since $\pi^n \in \sA_0$, we have $\pi^n_{\tilde t_n+} \le \f(\tilde t_n, \pi^n_{\tilde t_n}, q) \le y$,  and thus $\pi^n_{\tilde t_n+} = y$. Note that $\sM[v](\tilde t_n,X_{\tilde t_n}, y,q-y)>0$, then there exists $\e_n > 0$ such that $\sM[v](s, X_s,  \pi^n_s,q-\pi^n_s)>0$ for all $s\in (\tilde t_n, \tilde t_n + \e_n)$. This implies that $\pi^n_s \in O(s, X_s, q)$ and $\f(s, \pi^n_s, q) = \pi^n_s$. Now recall again that  $\pi^n \in \sA_0$, then we have $d (\pi^n)^c_s = 0$ and $\D \pi^n_s = 0$ for all $s\in (\tilde t_n, \tilde t_n + \e_n)$. Therefore, $\pi^n_s = y$ for all $s\in (\tilde t_n, \tilde t_n + \e_n)$, contradicting with the definition of $\tilde t_n$.

Summarizing, we have shown that every totally ordered subset of $\sA_{0}$ has an upper bound. Therefore, applying Zorn's Lemma, we 
conclude that $\sA_{0}$ has a maximal element in $\sA_0$, denoted  by $\pi^*$.
We claim that $\pi^*$ does satisfy \reff{suff}. Indeed, by its construction it suffices to prove 
\bea
\label{pi+eq}
\pi^*_{t+}=\phi(t,\pi^*_{t},q), \q \forall t\in[0,\hat\t_1), ~\mbox{$\hP$-a.s.}
\eea

Suppose not, then $c := \phi(t,\pi^*_{t},q) - \pi^*_{t+} > 0$. Define 
\beaa
\t &:=& \inf\{s>t: \pi^*_{s}\geq\pi^*_{t+}+c\}\wedge\hat\t_1,\\
\hat\pi^*_{s} &:=& \pi^*_s{\bf 1}_{[0, t]}(s) + [\pi^*_{t+}+c]{\bf 1}_{(t, \t]}(s)+[\pi^*_s\vee (\pi^*_{t+}+c)] {\bf 1}_{(\t, T]}(s), ~ s\in [0,T].
\eeaa
It is straightforward to check that $\hat\pi^* \in \sA_0$,  $\pi^* \prec \hat\pi^*$,  and $\pi^*_s < \hat\pi^*_s$ for $ s\in (t, \t]$.  This contradicts the fact that $\pi^*$ is a maximum element of $\sA_0$. This proves (\ref{pi+eq}), whence  the
proposition.
\qed

We are now ready to state the man result of this section.
\begin{thm}
\label{thm-existence}
Assume all the conditions of Proposition \ref{prop-verification} hold. Then $v=V$ and there exists an optimal strategy $\pi^* \in \sA(0,0,q)$ such that $v(0,x,0,q) = J^1(0,x,0,q; \pi^*)$.
\end{thm}

{\it Proof.} Combining Propositions \ref{prop-suff} and \ref{prop-existence}, there exists $\pi^*\in \sA(0,0,q)$ such that \reff{epi0} holds. Repeating the same arguments for each $n$, we may extend $\pi^*$ appropriately on $[0, \hat\t_n]$ such that 
\beaa
&&\sum_{i=0}^{n-1}  \mathbb{E}\Big\{\int_{\hat\t_i}^{\hat\t_{i+1}}  \sL[v](s,  X_s, \pi^i_s,   Q^{\pi,i}_s) ds + \int_{\hat\t_i}^{\hat\t_{i+1}} \sM[v] (s,  X_s,  \pi_s,  Q^{\pi}_s) d\pi^c_s \\
&& +\sum_{\hat\t_i\le s<\hat\t_{i+1}}  \int_0^{\D \pi_s} \sM[v](s,  X_s,  \pi_s + u, Q^{\pi}_s - u) du \Big\} = 0,
\eeaa
which, following the proof of Proposition \ref{prop-verification}, implies that
\beaa
\mathbb{E}\Big\{\int_0^{\hat\t_n}U(X_{s},Q^{\pi^*}_{s})d(\pi^*)^c_s+\sum_{0\leq s<\hat\t_n}D(X_s,Q^{\pi^*}_{s},\Delta\pi^*_{s})+ v(\hat\t_n, X_{\hat\t_n}, \pi^*_{\hat\t_n}, Q^{\pi^*}_{\hat\t^n})\Big\} =  v(0, X_0, 0, q).
\eeaa
Sending $n\to\infty$, and recalling the terminal condition in \reff{terminal}, we see that  
\beaa
v(0, X_0, 0, q) = J^1(0,x,0,q;\pi^*) \ge V(0,x,0,q).
\eeaa
This, together with Proposition \ref{prop-verification}, completes the proof.
\qed

\begin{rem}
{\rm
 Based on Proposition \ref{prop-suff}  we can roughly describe the optimal strategy $\pi^*$ as follows. 
At each time $t\in [\hat\t_i, \hat\t_{i+1}]$ between the two jump times of $N$,  there is an ``inaction region" 
 $O(t,X_{t}, Q^{\pi^*}_{\hat\t_i})$, which  
 is an open set, and therefore can be decomposed into open intervals. If $\pi^*_t - \pi^*_{\hat\t_i}\in O(t,X_{t}, Q^{\pi^*}_{\hat\t_i})$, then it stays ``flat." If it is
 at the boundary of  $O(t,X_{t},Q^{\pi^*}_{\hat\t_i})$, hence the boundary of one of the open intervals, then it either jumps to 
 $\f(t,X_t, Q^{\pi^*}_{\hat\t_i})$, i.e, the boundary of nearest neighboring interval above it, if 
 $\f(t,X_t, Q^{\pi^*}_{\hat\t_i})>\pi^*_t$, or move along with the boundary of $O(t,X_{t}, Q^{\pi^*}_{\hat\t_i})$, when
 $\f(t,X_t, Q^{\pi^*}_{\hat\t_i})=\pi^*_t$. In particular, when  $O(t,X_{t}, Q^{\pi^*}_{\hat\t_i})$ is simply connected, then
 $\pi^*$ essentially behaves like an optimal singular stochastic control. However, it is not clear to us that $O(t,X_{t}, Q^{\pi^*}_{\hat\t_i})$ will be simply connected, and consequently the optimal strategy may jump multiple (even infinitely many) times between $[\hat\t_i, \hat\t_{i+1}]$.
%
\qed}
\end{rem}


\begin{thebibliography}{1}

\bibitem{Alfonsis1}
Alfonsi, A., Fruth, A., Schied, A. (2008), {\it Constrained portfolio liquidation in a limit order book model}. 
Advances in mathematics of finance, 9Ð25, Banach Center Publ., {\bf 83}, Polish Acad. Sci. Inst. Math., Warsaw, 

\bibitem{Alfonsis2}
Alfonsi, A., Fruth, A., and Schied, A. (2010), \emph{Optimal execution
strategies in limit order books with general shape functions}. Quant. Finance {\bf 10},
no. 2, 143-157.

\bibitem{AlfonsiSchied}
Alfonsi, A. and Schied, A. (2010), \emph{Optimal execution and
absence of price manipulations in limit order book models}. SIAM J. Financial Math.
{\bf 1}, 490-522.

\bibitem{Alfonsis}
Alfonsi, A.,  Schied, A., and Slynko, A. (2012), \emph{Order Book
Resilience, Price Manipulation, and the Positive Portfolio Problem}.
SIAM J. Finan. Math., {\bf 3}, no. (1), 511-533.

\bibitem{Avellaneda}
Avellaneda, M. and Stoikov, S.  (2008), \emph{High-frequency trading in a limit order book}. Quantitative Finance, {\bf 8}, no. 3, 217-224.

\bibitem{BayLud}
Bayraktar, E. and Ludkovski, M. (2010), \emph{Optimal trade execution in illiquid markets}.
Mathematical Finance, 21: no. doi: 10.1111/j.1467-9965.2010.00446.x

\bibitem{Bieleckis2004}
Bielecki, T., Jeanblanc, M., and Rutkowski, M. (2004),
\emph{Hedging of Defaultable Claims}. Paris-Princeton Lectures on
Mathematical Finance 2003, 1-132, {\sl Lecture Notes in Math.}, {\bf 1847}, Springer, Berlin.


\bibitem{Cetins}
Cetin, U.,  Jarrow, R. A., and Protter, P. (2004), \emph{Liquidity
Risk and Arbitrage Pricing Theory}. Finance and Stochastics, {\bf 8}, no. 3, 311-341.


\bibitem{CIL} Crandall, M.G., Ishii, H., and Lions, P-L. (1992) {\it User's guide to
viscosity solutions of second order partial differential equations}, {\sl
Bull. Amer. Math. Soc. (NS)}, {\bf 27}, 1--67.

\bibitem{ETZ}
Ekren, I., Touzi, N., and Zhang, J.  {\it Viscosity Solutions of Fully Nonlinear Parabolic Path Dependent PDEs: Part I}, preprint,  arXiv:1210.0006.


\bibitem{Foucault}
Foucault, T., Kadan, O., and Kandel. E. (2005), \emph{Limit Order Book as a Market for Liquidity}. Review of Financial Studies, {\bf 18}(4), 1171-1217.

\bibitem{Gatherals}
Gatheral, J.,  Schied, A., and Slynko, A. (2012), \emph{Transient
linear price impact and Fredholm integral equations}. Mathematical Finance, {\bf 22}, no. 3, 445-474.



\bibitem{Handa}
Handa, P. and Schwartz, R. A. (1996), \emph{Limit Order Trading}. Journal of Finance, {\bf 51}, no. 5, 1835-1861.

\bibitem{Hollifield}
Hollifield, B., Miller, R. A.,  and Sandas, P. (2004),  \emph{Empirical Analysis of Limit Order Markets}. Review of Economic Studies, {\bf 71}, no. 4, 1027-1063.



\bibitem{OW}
Obizhaeva, A. and Wang, J. (2013), \emph{Optimal Trading Strategy and
Supply/Demand dynamics}. Journal of Financial Markets, {\bf 16}, no. 1, 1-32.


\bibitem{Rosu}
Rosu, I. (2009), \emph{A Dynamic Model of the Limit Order Book}. The
Review of Financial Studies,  {\bf 22}, 4601-4641.

\bibitem{shreve}
Predoiu, S., Shaikhet, G., and Shreve, S.E. (2011), {\it Optimal Execution in a General One-Sided Limit-Order Book}, SIAM J. Finan. Math.,  {\bf 2}, no. 1, 183-212.



\end{thebibliography}
\end{document}